\documentclass{fic-l}
\usepackage{amscd}
\usepackage[mathscr]{eucal}
\let\mathcal=\mathscr
\usepackage{xy}
\xyoption{all}

\newtheorem{thm}[subsection]{Theorem}
\newtheorem{cor}[subsection]{Corollary}
\newtheorem{conj}[subsection]{Conjecture}
\newtheorem*{thm*}{Theorem}
\newtheorem*{prop*}{Proposition}
\newtheorem*{conj*}{Conjecture}
\newtheorem*{lemma*}{Lemma}
\newtheorem*{cor*}{Corollary}
\newtheorem*{claim*}{Claim}
\theoremstyle{definition}
\theoremstyle{remark}
\newtheorem*{rem*}{Remark}
\newtheorem*{example*}{Example}
\numberwithin{equation}{section}
\newcommand{\op}{\operatorname}

\newcommand{\bh}{\begin{*hypo}}
\newcommand{\eh}{\end{*hypo}}
\newcommand{\bea}{\begin{eqnarray}}
\newcommand{\eea}{\end{eqnarray}}
\newcommand{\be}{\begin{equation}}
\newcommand{\ee}{\end{equation}}

\newcommand{\slp}{\partial\hspace{-1ex}/}

\newcommand{\cam}{{\mathcal M}}

\newcommand{\cag}{{\mathcal G}}

\def\sym{\mathop{\mathrm{Sym}}\nolimits}

\def\ind{\mathop{\mathrm{Ind}}\nolimits}
\def\ker{\mathop{\mathrm{Ker}}\nolimits}

\def\C{{\mathbb C}}
\def\Q{{\mathbb Q}}

\def\lc{\mathop{\mathcal{C}}\nolimits}

\def\Z{{\mathbb Z}}

\def\R{{\mathbb R}}

\title{Heegaard Splittings and Seiberg-Witten Monopoles}
\author{Yi-Jen Lee} 
\address{Department of Mathematics, Purdue University
\\W. Lafayette, IN 47907 U.S.A.\\ \texttt{\small yjlee@math.purdue.edu}}
\thanks{supported in part by NSF grant DMS \#0333163}

\date{%!! Preliminary version!!}
September 2004}
\begin{document}
\maketitle \thispagestyle{empty}
\begin{abstract}
This is an expansion on my talk at the Geometry and Topology
conference at McMaster University, May 2004. 

We outline a program to relate the Heegaard Floer homologies of
Ozsvath-Szabo, and Seiberg-Witten-Floer homologies as defined by 
Kronheimer-Mrowka. The center-piece of this program is the 
construction of an intermediate version of Floer theory, which exhibits
characteristics of both theories. 
\end{abstract}

\section{The conjecture}
This is a preliminary report on a long program aiming at a proof of
the conjectural equivalence between the Heegaard Floer homologies of 
Ozsvath-Szabo, and the monopole Floer homologies as defined by
Kronheimer-Mrowka.

Besides giving an overall picture, we give in \S5, 6, 8, and 9.1 a 
survey of some partial results towards this goal, with details
deferred to papers in preparation \cite{L:model, L:moduli}.

Throughout this article, we shall work with an unspecified coefficient
ring \(R\), which may be \(\Z/2\Z, \Z, \Q\), or \(\R\).
Take \(R=\Z/2\Z\) if desired, since we ignore 
the orientation issue in this article.  

Let \(Y\) be a compact oriented 3-manifold and \(\frak{s}\) be a
spin-c structure on \(Y\).
In \cite{OS1}, Ozsvath-Szabo defined four versions of Heegaard Floer
homologies associated to \((Y, \frak{s})\), \(\op{HF}^-(Y, {\frak
  s})\), \(\op{HF}^\infty(Y, {\frak s})\), \(\op{HF}^+(Y, {\frak s})\), and
\(\widehat{\op{HF}}(Y, {\frak s})\).

Let \(\op{HM}\) denote the Seiberg-Witten-Floer homologies defined by
Kronheimer-Mrowka \cite{KM:floer}. 
They come in three flavors, \(\hat{\op{HM}}\),
\(\bar{\op{HM}}\), \(\check{\op{HM}}\). Basing on Kronheimer-Mrowka's
construction, we shall introduce in \S 5 below 
a fourth version, \(\op{HM}^{tot}\), in parallel to
\(\widehat{\op{HF}}\). In addition to the pair \((Y, \frak{s})\),
these Seiberg-Witten-Floer homologies depend on the cohomology class
of perturbation two form, \([\omega]\), and we denote them by 
\(\op{HM}(Y, \frak{s}; [\omega])\).

In spite of their very different origins, these Floer homologies have
identical formal properties. They are both \(R[U]\)-modules, where
\(U\) is a degree \(-2\) chain map, and the first three flavors of
both Floer homologies fit into long exact sequences, which we call
the {\em fundamental exact sequences}:
\begin{gather*}
\cdots \to \op{HF}^-\to \op{HF}^\infty\to \op{HF}^+\to \cdots\\
\cdots\to \hat{\op{HM}}\to 
\bar{\op{HM}}\to \check{\op{HM}}\to \cdots
\end{gather*}

\begin{conj}
Let \((Y, \frak{s})\) be as the above, and let 
\begin{equation}\label{w-class}
[w]=2\pi c_1({\frak
  s}).
\end{equation} 
Then there are isomorphisms of \(R[U]\)-modules
\begin{gather*}
\op{HF}^-(Y, {\frak s})\simeq\hat{\op{HM}}(Y, {\frak s}; [w]),\\
\op{HF}^\infty(Y, {\frak s})\simeq\bar{\op{HM}}(Y, {\frak s}; [w]),\\
\op{HF}^+(Y, {\frak s})\simeq\check{\op{HM}}(Y, {\frak s}; [w]),\\
\widehat{\op{HF}}(Y, {\frak s})\simeq\op{HM}^{tot}(Y, {\frak s}; [w]),
\end{gather*}
which are natural with respect to the fundamental exact sequences of 
Heegaard and Seiberg-Witten Floer homologies. 
\end{conj}

This conjecture has been verified for all known computations of both
sides. In addition, since both Heegaard Floer homologies and
Seiberg-Witten Floer homologies satisfy surgery exact sequences, if
there is a map between two theories natural with respect to the
surgery exact sequence, then the conjecture holds. 

However, the difficulty in proving the above conjecture lies
precisely in finding such a natural map. A quick look at the
construction of both theories finds them very different both in
geometric contents and abstract frameworks. We shall return to 
this subject in \S 3.

There are many variants and extensions of this conjecture
which we omit in this article. For example, twisted versions of
both Floer homologies are conjectured to relate in a similar fashion. In
addition, Floer homologies are the building block for the
definition of 4-manifold invariants and contact invariants in both
theories, and these invariants are also conjectured to be equal. 
Overall, the Seiberg-Witten theory is more closely related to the
geometry of underlying manifolds (e.g. scalar curvature), while the
Ozsvath-Szabo theory, more combinatorial in flavor, is in general more
computable.

\subsection{Basic ingredients of Seiberg-Witten-Floer theory}

Due to limitation of space, we shall not explain the construction of
either theories, but refer the reader to the original literature. 
Here we shall only recall some basic notions for the sake of fixing
notation and terminology.

A Seiberg-Witten {\em configuration} is a pair \((A, \psi)\), where
\(\psi\) is a section of the spinor bundle \(S\) over 3- or
4-dimensional spin-c manifold, and \(A\) is a connection on \(\det S\).
Because both the 3-dimensional and the 4-dimensional
contexts appear in this article, we shall reserve the unhatted 
notation \((A, \psi)\) for
3-dimensional configurations, and the hatted version \((\hat{A},
\hat{\psi})\) for 4-dimensional Seiberg-Witten configurations. 
  
By a {\em Seiberg-Witten-Floer theory on \((Y,
\frak{s})\) perturbed by \(\omega\)}, we mean the following. 
As a formal Morse theory, its chain groups are generated by what we
call ``{\em Seiberg-Witten critical points}'', which are
(gauge-equivalence classes of) solutions to 
the 3-dimensional (perturbed) Seiberg-Witten equations: 
\[
\left\{ \begin{array}{l}
\slp_A \psi =0+\cdots \\
F_A-i\sigma (\psi, \psi)=i\omega+\cdots, \end{array}\right.
\hskip 6cm
(SW_3)\]
where \(\sigma(\psi, \psi)\) is certain quadratic function of
\(\psi\), and \(\omega\) is a closed 2-form.

The boundary map of the Floer complex is defined by counting what we
call ``{\em Seiberg-Witten flow lines}'' between two Seiberg-Witten critical
points, which are 
(gauge-equivalence classes of) solutions to 
the 4-dimensional (perturbed) Seiberg-Witten equations on \((\R\times
Y, \frak{s})\): 
\[
\left\{ \begin{array}{l}
\slp_{\hat{A}} \hat{\psi }=0+\cdots \\
\op{SD}(F_{\hat{A}})-i\sigma (\hat{\psi}, \hat{\psi})=i\op{SD}(\omega)+\cdots, \end{array}\right.
\]
where \(\op{SD}\) denotes the self-dual part of a
2-form, and \(\omega\) means the pull-back of the two form \(\omega\)
on \(Y\) to \(\R\times Y\), and \(\frak{s}\) now denotes a spin-c
structure on \(\R\times Y\) via the identification of the spaces of
spin-c structures on \(Y\) and \(\R\times Y\).

The dots in the Seiberg-Witten equations above indicate that
additional perturbation is needed to achieve transversality. These
perturbations %are possibly non-local, but 
will not change the
cohomology class of the CSD functional (or, in the terminology of \S3,
the class of the homomorphism \(\op{PR}\)). This technical point will
be omitted in this article; for the precise form of these
perturbations, see \cite{KM:floer}.

The 3-dimensional Seiberg-Witten invariant, denoted \(\op{Sw}_3\), is
the Euler characteristic of the Seiberg-Witten-Floer homology. When 
\(b_1(Y)\geq1\), it may be obtained by a straightforward signed count
of solutions to (\(SW_3\)). It is independent of the perturbation when
\(b_1(Y)>1\).  When \(b_1=1\), it depends on the chamber \([w]\) is
in. In this article, \(\op{Sw}_3\) will always mean the invariant in
the ``Taubes chamber'' .

\subsection{Basic ingredients of Heegaard Floer theory}

Let 
\[
Y=H_-\cup_{\Sigma} H_+, 
\]
be a Heegaard splitting of \(Y\), namely the 3-manifold \(Y\) is
separated into the two handlebodies \(H_+, H_-\) by a Heegaard 
surface \(\Sigma\) of genus \(g\).

A Morse function \(f: Y\to \R\) is said to {\em adapt to the Heegaard
  splitting} if \(f^{-1}(0)=\Sigma\) is the Heegaard surface, and
\[H_+=f^{-1}\R_{\geq0}, \quad H_-=f^{-1}\R_{\leq0}\] 
contain respectively one minimum \(p_+\) and \(g\) index 1 critical points,
and one maximum \(p_-\) and \(g\) index 2 critical points. 

Let \(\alpha_i, \, i=1, \ldots, g\) denote the descending cycles on
\(\Sigma \) from the \(g\) index 2 critical points. Similarly, let
\(\beta_i, \, i=1, \ldots, g\) denote the ascending cycles on
\(\Sigma \) from the \(g\) index 1 critical points. Let
\[{\mathbb T}_\alpha:=\alpha_1\times \cdots\times \alpha_g, \quad
{\mathbb T}_\beta:= \beta_1\times \cdots\times \beta_g \subset
\sym^g\Sigma.\]
Suppose \({\mathbb T}_\alpha\) and \({\mathbb T}_\beta\) intersects
transversely. 

As a formal Morse theory, the chain groups in the Heegaard Floer
theory is generated by what we call the ``{\em Heegaard critical
  points}'', which are intersection points of \({\mathbb T}_\alpha\) and
\({\mathbb T}_\beta\). The boundary map is defined by counting what we
call ``{\em Heegaard flow lines}'', which are holomorphic disks
\[
\mu: \R\times [0,1]\to \sym^g\Sigma, \quad \text{with \(\mu(\cdot,0)\in
  {\mathbb T}_\alpha\), \(\mu(\cdot,1)\in
  {\mathbb T}_\beta\)}.
\]
Of paramount importance in Heegaard Floer theory is the choice
of a {\em base point}: 
Let \(z\in \Sigma\) be a point avoiding the descending and ascending
cycles, and let \(\gamma_z\subset Y\) be the flow
line of \(f\) from \(p_+\) to \(p_-\) through \(z\). We shall explain
the role of \(z\) and \(\gamma_z\) in Heegaard
Floer homologies in \S4. 

\section{First motivation: Taubes's work on \(\op{Sw}=\op{Gr}\)}

In fact, it is not surprising that the Seiberg-Witten-Floer homologies
should be related to curve-counting invariants. Since Taubes's seminal
work, the relation between Seiberg-Witten theory and Gromov's theory
of pseudo-holomorphic curves has been well-known.

\begin{thm}[Taubes]
Let \((X, \varpi)\) be a closed, oriented, symplectic 4-manifold, and
\(\frak{s}\) be a spin-c structure on \(X\). Then 
\[
\op{Sw}_4(X, \frak{s})=\op{Gr}(X, \frak{s}),
\]
where \(\op{Sw}_4\) is the Seiberg-Witten invariant for 4-manifolds,
and \(\op{Gr}\) is a variant of Gromov invariant that counts embedded,
possibly disconnected, pseudo-holomorphic curves (with multiplicity) 
in the homology class determined by \(\frak{s}\). 
\end{thm}

In the case where \(X\) is an algebraic manifold, this is just a
simple analog of the correspondence between line bundles and divisors,
and has been known since the first discovery of Seiberg-Witten theory.
We shall briefly explain some of Taubes's ideas, as it will be
central to our program. 

First, choose a metric on \(X\) with respect to which the symplectic
form \(\varpi\) is self dual (hence harmonic).
The metric, together with \(\varpi\), determines an almost complex
structure on \(X\).
Observe that the Clifford action by \(\varpi\) splits the spinor
bundle into a direct sum of eigenspaces: 
\begin{equation}\label{S-split}
S=E\oplus E\otimes K^{-1},
\end{equation}
where \(K^{-1}\) is the anti-canonical bundle.
We shall therefore write
\[
\hat{\psi}=(\hat{\alpha}, \hat{\beta})
\]
in accordance with this splitting. As the Levi-Civita connection
determines a connection of \(K^{-1}\), specifying a spin-c connection
is equivalent to specifying a connection \(\hat{A}^E\) of \(E\). We shall therefore
denote a Seiberg-Witten configuration in this context as
\[
(\hat{A}^E, (\hat{\alpha}, \hat{\beta})).
\]

Taubes considered perturbation to the Seiberg-Witten equations on
\(X\) by the two-form 
\[
\omega=r\varpi+\text{insignificant terms,} \quad \text{where \(r\gg1\)
  is a constant}.
\]
He showed that as \(r\to \infty\), the zero locus \(\alpha^{-1}(0)\)
approaches a holomorphic curve, which is also the support of the 
current \(\lim_{r\to \infty} F_{A^E}\). 
We shall call this a ``{\em Taubes curve}''.

His proof also gives a local description of the large \(r\)
Seiberg-Witten solutions: Away from the Taubes curve, the solutions
approximate the simple form
\[
(\hat{A}^E, (\hat{\alpha}, \hat{\beta}))=(0, (\sqrt{r}, 0)), \quad
\text{with respect to a trivialization of \(E\),}
\]
while near the Taubes curve, it approximates a family of vortex
solutions parameterized by the Taubes curve. 

A heuristic way to understand this result is as follows. Locally, 
varying \(r\) has the effect of rescaling by \(\sqrt{r}\). Thus, as
\(r\to \infty\), the Nijenhuis tensor is ignorable, reducing to the simpler
K\"{a}hler situation mentioned above. 

As we shall need several variants of Taubes's settings, we shall refer
to them in general as ``{\em Seiberg-Witten-Taubes}'' theories,
abbreviated as SWT. 

\subsection{Taubes's picture in 3 dimensions: Morse-Novikov theory}

\noindent\ \par\bigbreak
\noindent{\bf (i) SWT theory on closed 4-manifolds.}
Instead of a symplectic form, one may consider any self-dual harmonic
2-form for \(\varpi\). For generic metric, such \(\varpi\) vanishes
along a set of circles in \(X\), and \(K^{-1}\) and the splitting
(\ref{S-split}) makes sense away from these circles. Taubes proved in \cite{T2}
an extension of the convergence theorem in part I of \cite{T} to this 
setting. Here, a Taubes curve is a pseudo-holomorphic subvariety with
boundary at \(\varpi^{-1}(0)\); more precisely, the intersection
number of the Taubes curve with any linking 2-sphere of \(\varpi^{-1}(0)\) is 1.

However, lack of understanding on the local behaviors of 
the Seiberg-Witten solutions and the corresponding pseudo-holomorphic 
subvarieties deters progress on establishing a full equivalence
generalizing Theorem 2.1. In fact, even the desired generalization of 
the Gromov invariant is undefined. 
\bigbreak 

\noindent {\bf (ii) SWT theory on closed 3-manifolds.}
The 3-dimensional story is considerably simpler. Consider perturbation
to the 3-dimensional Seiberg-Witten equations on a closed spin-c
3-manifold \((Y, \frak{s})\) by \(r\varpi+\cdots\), where \(\varpi\)
is now a harmonic 2-form. The Hodge dual \(*\varpi=:\vartheta\) is a
harmonic 1-form; for generic metric, it is Morse. 
This brings us to the realm of Morse-Novikov theory, as follows.

The zero locus \(\varpi^{-1}(0)\) is now the critical points of
\(\vartheta\): they come in \(g\) pairs of index 2 and index 1
critical points. Let \(\check{\vartheta}\) denote the dual vector field of
\(\vartheta\). The counterpart of a Taubes curve is a union of
finite-length flow lines of \(-\check{\vartheta}\), such that the union
of their boundary is precisely the critical set \(\varpi^{-1}(0)\).
We call this a ``{\em Taubes orbit}'', and the flow lines constituting
a Taubes orbit the ``{\em constituent flow lines}''.
%\footnote{many of the references below should be changed accordingly}

Basing on Taubes's picture, we wrote down a counting invariant of
Taubes orbits, \(I_3\), and conjectured:
\begin{conj*} {\rm \cite{HL1}}
Let \((Y, \frak{s})\) be a closed spin-c 3-manifold with
\(b_1(Y)>0\). Then \[\op{Sw}_3(Y, \frak{s})=I_3 (Y, \frak{s}).\]
\end{conj*}
In fact, \(I_3\) is just a special case of an invariant \(I\) in general
Morse-Novikov theory. This invariant \(I\) takes the form of a product
of the Reidemeister torsion of the Morse-Novikov complex, with a
dynamical zeta function that counts closed orbits of the flow of
\(-\check{\vartheta}\). Moreover, we showed:
\begin{thm*}{\rm \cite{HL1, HL2}}
For an oriented, closed manifold \(M\), \[I(M)=\tau(M),\]
\(\tau(M)\) being the combinatorially defined Reidemeister torsion of
the manifold \(M\). In particular, 
\(I_3 (Y)=T(Y)\), the Turaev torsion. 
\end{thm*}
It is known that \(\op{Sw}_3\) is equivalent to the Turaev torsion,
either by the surgery formulae of Meng-Taubes, or by the TQFT
arguments of Donaldson-Mark. Thus, the above Conjecture is proven
indirectly. 
\bigbreak

\noindent{\bf (iii) SWT Floer theory on closed 3-manifolds.}
It is also possible to generalize this picture to Floer theory. 
Consider the Seiberg-Witten-Floer theory with the same perturbation
\(\omega\). Its chain group will be generated by 3-dimensional
SWT solutions, which corresponds to the Taubes orbits 
described above. The boundary map is defined
via the moduli spaces of Seiberg-Witten solutions on \(\R\times Y\), with
perturbation \(\op{SD}(\omega)+\cdots\). The latter now vanishes along lines
\(\R\times \varpi^{-1}(0)\). The corresponding Taubes curves are now
%Taubes's picture predicts that these
%Seiberg-Witten solutions correspond to Taubes curves, which are now
pseudo-holomorphic curves in \(\R\times Y\backslash (\R\times
\varpi^{-1}(0))\), with boundary at the vanishing locus \(\R\times \varpi^{-1}(0)\).
Thus, the analog of Theorem 2.1 states that the Seiberg-Witten-Floer
homology is equivalent to a symplectic version of Floer homology,
analogous to the contact homology.

In the simplest special case, when \(Y\) is mapping tori, it is
possible to choose \(\vartheta\) such that
\(\varpi^{-1}(0)=\emptyset\). Hutchings established some foundations of
the proposed symplectic version of Floer homology in this case, which
he calls ``periodic Floer homology''.

While this picture 
potentially offers an interesting connection with contact
homology, it has so far not offered great help on the understanding of
Seiberg-Witten-Floer homologies, as its construction and computation
is no simpler than Seiberg-Witten-Floer homologies. 
 
\subsection{Heegaard-splittings and real-valued Morse theory}
Instead of the Morse-Novikov theory of closed 1-forms, 
it is desirable to have a
geometric interpretation of 3-dimensional Seiberg-Witten theory in
terms of real-valued Morse theory, for the following two reasons:
\begin{itemize}
\item  The important case of rational homology spheres (\(b_1(Y)=0\))
  is excluded from the discussion in \S2.2.
\item  Via the Heegaard splitting associated to a real-valued Morse
  function, it is possible to reduce the computation to 2-dimensions,
  namely, to the Heegaard surface. This would entail great simplifications.
\end{itemize}
At the first sight, this does not appear possible.
Suppose \(\vartheta=df\) for a real valued function \(f\). The
harmonicity of \(\vartheta\) implies the harmonicity of \(f\), but 
there is no non-constant harmonic function on 
a closed manifold.

This, however, does not constitute a serious obstacle. 
Real-valued harmonic functions
do exist on {\em non-compact} 3-manifolds.  For instance, one may
consider
\begin{enumerate}
\item  {\bf (3-d Euclidean SWT theory)} \(Y\# \R^3\) with Euclidean metric at
  infinity. Seiberg-Witten-Taubes theories on such 3-manifolds are
  considered in \cite{L1}.
\item {\bf (3-d cylindrical SWT theory)}
Deleting two points \(p_+, p_-\)
  from \(Y\), then choosing a complete metric on the resulting open
  manifold, so that it has two cylindrical ends \(\R_{\pm}\times
  S^2\). There exist harmonic functions \(f\) on such cylindrical
  manifolds which are asymptotic to 
\begin{equation}\label{asym-tau}
\tau+C_\pm\end{equation} 
on the two
  cylindrical ends, where \(\tau \) parameterizes \(\R_\pm\). Consider
  spin-c structures such that at the spheres at infinity,
  \(E|_{S^2_{\pm\infty}}\) has nonnegative degree. In particular, when
  the degree is 0, denote the corresponding spin-c structure by
  \(\frak{s}^\circ\). More generally, one may consider any 3-manifolds
  with a finite number of cylindrical ends, with a harmonic function
  with asymptotic condition (\ref{asym-tau}) on the ``negative'' and
  ``positive'' ends as \(\tau\to-\infty\), \(\infty\) respectively. 
\end{enumerate}
Taubes's pictures for both settings are similar; so we shall take the
second situation for example. Let \(f\) be a Morse function adapted to
a Heegaard splitting, as described in \S1.3. 
Then the constituent flow lines of a Taubes orbit of 
spin-c structure \(\frak{s}^\circ\)
may be described as intersection points of descending
cycles \(\alpha_i\) and ascending cycles \(\beta_j\), and the Taubes orbit
corresponds precisely to a intersection point of \({\mathbb
  T}_\alpha\) and \({\mathbb T}_\beta\) in \(\sym^g\Sigma\). Notice
that the fact that \(f\) is real-valued excludes the possibility of
closed orbits from constituent flow lines; in particular, 
all constituent flow lines in Taubes orbits has multiplicity 1. 
This represents another
advantage of real-valued Morse theory over the Morse-Novikov theory. 

Moreover, as explained in \S2.2, the Taubes curves corresponding to
SWT flow lines now correspond to holomorphic curves
with boundary along the lines \(\R\times \op{Crit}(f)\).
Analogy with the Atiyah-Floer conjecture (cf. \S7.2, 9.2) leads to the
expectation that these corresponds to holomorphic disks in
\(\sym^g\Sigma\) with boundary along \({\mathbb T}_\alpha\),
\({\mathbb T}_\beta\), and that the Seiberg-Witten-Floer homology 
corresponds to the Floer homology of ``Lagrangian'' intersections of 
\((\sym^g\Sigma; {\mathbb T}_\alpha, {\mathbb T}_\beta)\).
This is precisely what Ozsvath-Szabo calls ``\(\op{HF}^\infty\)''.

Independent of its relation with Seiberg-Witten theory, this
Morse-theoretic picture defines certain topological invariants. However,
more serious problems of these ideas are:
\begin{itemize}
\item These are decidedly different from the ordinary Seiberg-Witten 
invariant or Floer homologies of \(Y\).
\item In fact, they are not interesting invariants, since they only
  depends on homological information of \(Y\).
\end{itemize} 
%The Seiberg-Witten invariants defined from these non-compact manifolds
%may or may not agree with the Seiberg-Witten invariants of \(Y\);
%they nevertheless promise to define topological invariants of \(Y\).

The first point is not unexpected, since this story supposedly
corresponds to Seiberg-Witten theory on {\em noncompact} manifolds,
while gauge-theoretic invariants on noncompact manifolds are typically
very sensitive to asymptotic conditions. In fact, viewing the
cylindrical manifolds in situation 2 as the closed manifold \(Y\)
minus two points, the choice of spin-c structures and perturbations in
the Seiberg-Witten-Taubes theory are those which do not extend across the two
points.

The second problem is, of course, a lot worse. To illustrate it,
notice that 
\[\begin{split}
\# & \,  ({\mathbb T}_\alpha \cap {\mathbb T}_\beta)\\
& =\det \, (\# (\alpha_i\cap \beta_j))\\
& =\begin{cases}0 &\text{as \(b_1(Y)>0\)};\\ |H_1(Y; \Z)| &\text{as
      \(b_1(Y)=0\),}
\end{cases}
\end{split}\]
which contains a lot less information than the ordinary Seiberg-Witten
invariants of 3-manifolds.

This is where people abandoned this approach. 
It took the advent of Ozsvath-Szabo's amazing idea--filtration-- to revive this
program. (See \S4).

\section{The Floer-theory framework}
We now examine the abstract frameworks of both Floer theories. 

\subsection{Basic ingredients of a Floer theory}
A typical Floer theory models on  a formal Morse theory over
an infinite dimensional space \(\lc\). Usually, one has
a possibly non-exact Morse 1-form over \(\lc\). Furthermore, different
from the finite-dimensional Morse theory, the (relative)
index of the critical points is only defined modulo an integer \(N\).
There is a minimal abelian covering \(\tilde{\lc}\to \lc\) such that
the Morse 1-form lifts
to a differential of a \(\R\)-valued Morse function, and indices of lifts of
the critical points are well-defined in \(\Z\) {\em for all parameters
of the theory} (e.g perturbations, metrics, spin-c
structures). Let \(G\) be the covering group, and let  
% corr to subgroups \(N\subset\ker \op{PR}\): covering group \(G/N\)
\[
\op{PR}: G\to \R; \quad \op{SF}: G\to \Z
\]
denote homomorphisms defined by the change in the values
of the Morse function and the indices under deck transformation.
(\(\op{PR}\) stands for ``period'', and \(\op{SF}\) stands for
``spectral flow''). 

These are important basic ingredients of the Floer
theory, because Floer complexes with
different local coefficients model on Morse complex on sub-covers
\(\hat{\lc}\) of \(\tilde{\lc}\), on which the Morse 1-form lifts to
be exact.
Thus, \(G, \op{SF}, \op{PR}\) encode the module structure,
coefficient rings, and grading of these various Floer complexes.
%\indent\(\leftrightarrow\) subgroups \(N\) of \(\ker \op{PR}\) \newline
%\indent \(\leftrightarrow\) 
%\remfont{(finitely generated \(R[G/N]\) modules, coefficients: \(R[\ker
%  \op{SF}|_N]\), graded mod \(\op{gcd}(\op{SF}|_N)\))}

Unless otherwise specified, a ``Floer complex'' means the version
corresponding to a minimal \(\hat{\lc}\). %\(\leftrightarrow\) Floer complex with
It has ordinary \(R\)-coefficients when:
\begin{equation}\label{mono}
\op{PR}|_{\ker\op{SF}}=\emptyset.
\end{equation}
The finite dimensional Morse-Novikov theory also leads one to expect
the Floer homology to be invariant under variation of metric and
perturbations, as long as 
\begin{equation}\label{ray}
\text{\(\op{PR}\) remains on a positive ray from 0.}
\end{equation}

\subsection{Comparing the frameworks of two theories}
The following table compares the formal framework of 
the Seiberg-Witten and Heegaard Floer
theories in the above notation. It gives a first indication of the
many fundamental differences between the two theories.
\medbreak
{\small
\begin{center}
\begin{tabular}{|c||c|c|}
\hline
& Seiberg-Witten & Heegaard\\
\hline\hline
\(\lc\) & \(\bigl\{(A, \psi)\bigr\}/\cag\): quotient of the   &\(\Omega(\sym^g(\Sigma); {\mathbb T}_\alpha,
{\mathbb T}_\beta)\): space of  \\
 & spaces of configurations  & paths \(\gamma: [0,1]\to \sym^g(\Sigma)\), 
\\
&by gauge group action & with \(\gamma(0)\in {\mathbb T}_\alpha\), \(\gamma(1)\in {\mathbb
  T}_\beta\)\\
\hline
& \(\lc\) singular at fixed points (reducibles) & \(\lc\) smooth\\
\hline
%\(G\) & \(H^1(\lc; \Z)=\) & \(\pi_2({\bf x}, {\bf x})=\)\\
% & \(H^1(Y; \Z)\) & \(\Z\oplus H^1(Y; \Z)\) {\scriptsize (*)}\\
\(G\) & \(H^1(\lc; \Z)= H^1(Y; \Z)\) & \(\pi_2({\bf x}, {\bf x})=\Z\oplus H^1(Y; \Z)\) {\scriptsize (*)}\\
\hline
\(\op{PR}\) & \(2\pi (2\pi c_1(\frak{s})-[w])\) {\scriptsize (\dag)} &
\(-C \oplus 0\) {\scriptsize (\ddag)}\\
&& \(C\) is a positive constant\\
\hline
\(\op{SF}\) &\(c_1(\frak{s})\) &\(2\oplus c_1(\frak{s})\)\\
\hline
\end{tabular}
\end{center}
}
\bigbreak

\begin{rem*}
(*) Projection to the \(\Z\) factor is given by the intersection number of
the base point \(z\) with the
2-chain in \(\Sigma\) associated to the 1-cycle in
\(\Omega(\sym^g(\Sigma); {\mathbb T}_\alpha, {\mathbb T}_\beta)\).  

(\dag) Here and below, we regard a cohomology class in \(H^2(Y)\) as a
homomorphism \(H^1(Y)\to \R\) via Poincar\'{e} duality.

(\ddag) This assumes that \({\mathbb T}_\alpha\), \({\mathbb T}_\beta\) are
Lagrangian. However, as defined by Ozsvath-Szabo, they are typically
only totally real. This technical point is ignored here, as this section
is for motivational purposes only. 
\end{rem*}

As the covering group \(G\) of the Heegaard Floer theory contains an
extra \(\Z\)- factor than that of the Seiberg-Witten theory, to put
both theories on the same footing, one regards the Heegaard Floer
complex (\(\op{CF}^\infty\)) as an {\em infinitely generated}
\(R[H^1(Y; \Z)]\) module. In other words, instead of a Morse theory on
\(\lc\), model the Heegaard Floer
theory on the Morse theory on the infinite-cyclic covering of \(\lc\)
corresponding to the projection of \(G\) to \(\Z\).
We denote this infinite-cyclic covering by \(\hat{\lc}_z\), since the
projection is determined by \(z\).

On the other hand, notice that according to the table, 
if one chooses the class of perturbation form to be
\([w]=2\pi c_1(\frak{s})\), then the homomorphism \(\op{PR}\) in
Seiberg-Witten theory is 0. In this case, the homomorphisms
\(\op{PR}\) and \(\op{SF}\) in Seiberg-Witten theory agree with the 
restriction of their Heegaard counterparts to the \(H^1(Y;\Z)\)-factor.
This is a first motivation for the choice
(\ref{w-class}) in Conjecture 1.1. Furthermore, with this choice,
the condition (\ref{mono}) holds in both theories, and thus both Floer
complexes are of \(R\)-coefficients. 

With these explained, one may continue to observe the numerous formal
parallelisms of both theories. However, these parallel aspects 
come from entirely different origins.
In fact, both Floer theories require refinements of the basic Floer
theory framework outlined in \S3.1 above, but the main difference is
that the refinement
needed for each theory is based on different principles: the Heegaard
Floer theory relies on a filtration of the complex associated to the
infinite-cyclic covering, while the Seiberg-Witten theory is an 
\(S^1\)-equivariant theory. 

We shall discuss these different
refinements separately in the next two sections. For now, 
we continue the table of comparison that highlights the difference of
the refinements. 

\medbreak
{\small 
\begin{center}
\begin{tabular}{|c||c|c|}
\hline
Formal analogies & Seiberg-Witten & Heegaard\\
\hline
\hline
both are  & because this is an  &
because this models on \\
\(R[u]\)-modules & equivariant theory:  & Morse
theory on a \(\Z\)-cover: \\
&  \(U\) generates \(H^*(BS^1)\) & \(U\) generates deck transformation\\
\hline 
both complexes  are & because of  & because there are  \\ 
\(\infty\)-generated  & reducible critical points  & \(\infty\)-many critical points \\
\(R\)-modules & (fixed
points) & on a \(\Z\)-cover\\
\hline
long exact sequences & this is a & this is the \\
relate \(\hat{\op{HM}}\), \(\bar{\op{HM}}\), \(\check{\op{HM}}\), &
fundamental sequence &relative sequence\\
and \(\op{HF}^-\), \(\op{HF}^\infty\), \(\op{HF}^+\) & of
\(S^1\)-equivariant theory & associated with a filtration\\
\hline
both have & \(\op{HM}^{tot}\) models on  & \(\widehat{\op{HF}}\) models
on homology \\
a 4th version & homology of an \(S^1\)-bundle & of a
fundamental domain\\
\hline
\end{tabular}
\end{center}
}
\bigbreak

\section{Filtration in Heegaard Floer theory}
The key point in Ozsvath-Szabo's construction is a filtration argument
which is less commonly seen in the Morse/Floer theory literature. We shall 
put their construction in abstract formulation, since 
the same construction is needed again in \S6.

\subsection{1-cocycles, local coefficients, and infinite-cyclic coverings}
Let \(M\) be a manifold, and \(Z\) be a 1-cocycle in \(M\), or
equivalently, a codimension 1
cycle in \(M\). Suppose the cohomology class \([Z]\) is primitive; in
particular, nontrivial. Such \(Z\) defines a local system, \(\Gamma(Z)\), over \(M\), by
assigning each 1-chain in \(M\) its intersection number with \(Z\).
Let \(\tilde{M}_Z\) be the infinite-cyclic covering associated to
\(Z\), namely, the covering obtained by cutting \(M\) open along
\(Z\), and gluing \(\Z\) copies of such. The homology of \(M\) with
local coefficients \(\Gamma(Z)\) may alternatively be regarded as the
homology of the covering \(\tilde{M}_Z\). This is a graded
\(R[\Z]=R[t, t^{-1}]\) module. 

In the Morse-Novikov context, given a closed Morse 1-form \(\theta\)
on \(M\) with cohomology class \begin{equation}\label{Z-coh}
[\theta]=\alpha [Z] \quad \text{for
  \(\alpha\in \R\)},\end{equation}
one may define a Morse complex with local
coefficients, \(M_*(\theta, \Gamma(Z))\).  
This is a module of the completed ring \(R[t^{-1}, t]]\) when
\(\alpha>0\), of completion in the opposite direction, \(R[t,
t^{-1}]]\) when \(\alpha<0\). Here and below, \(t\) denotes a {\em
negative} generator of the deck transformation, in the sense that \(t\)
decreases the value of \(\tilde{f}\), where \(d\tilde{f}\) is the lift
of \(\theta\) to \(\tilde{M}_Z\).

\subsection{Filtration associated to semi-positive 1-cocycles}
We say that \(Z\) is
{\em semi-positive} with respect to \(\theta\) if 
the cohomology class \([Z]\) satisfies (\ref{Z-coh}), and \(Z(\gamma)\geq 0\) if
the path \(\gamma\) is a flow line of \(-\check{\theta}\).

The positivity condition implies that \(Z\) gives a
filtration on the Morse complex, and we may consider the associated
subcomplexes (as \(R\)-modules), and quotient complexes,
which are now \(R[t]\)-modules:

Let \(\tilde{Z}\) be a lift of \(Z\) in the infinite-cyclic cover
\(\tilde{M}_Z\). It divides \(\tilde{M}_Z\) into two halves. Let
\(\tilde{M}_Z^-\subset \tilde{M}_Z\) be the lower half (with respect
to the flow of \(-\check{\theta}\)).
%, \(\check{\theta}\) being the vector field dual to \(\theta\).
The subcomplex \(M^-_*(\theta, \Gamma(Z))\) may be understood as the
Morse complex for \(\tilde{M}_Z^-\). The quotient complex
\[
M^+_*(\theta, \Gamma(Z)):=M_*(\theta, \Gamma(Z))/M^-_*(\theta, \Gamma(Z))
\]
is then the complex of the pair \((\tilde{M}_Z, \tilde{M}_Z^-)\).

From the short exact sequence of \(R\)-modules
\[
0\to M^-_*(\theta, \Gamma(Z))\to M_*(\theta, \Gamma(Z))\to
M^+_*(\theta, \Gamma(Z))\to 0,
\]
we have a long exact sequence of the pair \((\tilde{M}_Z,
\tilde{M}_Z^-)\).

Notice that a different, though cohomologous choice of \(Z\) yields 
an equivalent local system, but a possibly different filtration. 
Thus, another semi-positive 1-cocycle \(Z'\) cohomologous to \(Z\) 
gives the same \(H_*(M_*(\theta, \Gamma(Z)))\), but possibly different
\(H_*(M_*^-(\theta, \Gamma(Z)))\) and \(H_*(M_*^+(\theta, \Gamma(Z)))\).
Similarly, \(H_*(M_*(\theta, \Gamma(Z)))\) is invariant under any 
exact perturbation to \(\theta\), but \(H_*(M_*^-(\theta,
\Gamma(Z)))\) and \(H_*(M_*^+(\theta, \Gamma(Z)))\) in general are only
invariant under {\em small} exact perturbations of \(\theta\).

\begin{example*}
Write \(\theta=df\) for a circle-valued function \(f\), and
take \(Z\) to be the 1-cocycle \(\theta\), or the codimensional 1
cycle given by a level set of \(-f\). Then \(Z\) is a semi-positive, 
and the associated filtration is the filtration by energy.
\end{example*}
\begin{example*}
Let \(P\) be a K\"{a}hler manifold, and consider 
a symplectic version of Floer theory on \(P\), \(\lc\) is a (relative) loop
space \(\Omega\) of \(P\). A 1-cycle in
\(\Omega\) traces out a (relative) 2-cycle in \(P\), and a flow line in
\(\Omega\) corresponds to a holomorphic curve in \(P\). The
intersection with a (complex) hypersurface in \(P\) thus defines a
semi-positive 1-cocycle. See \cite{seidel:icm}.
\end{example*}
\begin{example*}[Heegaard Floer homologies]
The Heegaard Floer theory is a variant of the construction of the
previous example. The hypersurface in this context is \(\{z\}\times \sym^{g-1}\Sigma\).
The first three versions of Heegaard Floer homologies are
respectively formal analogs of \(H_*(M^-_*(\theta, \Gamma(Z)))\),
\(H_*(M_*(\theta, \Gamma(Z)))\), \(H_*(M^+_*(\theta, \Gamma(Z)))\)
above, and the fundamental exact sequence of Heegaard Floer homologies
models on the relative exact sequence of the pair \((\tilde{M}_Z,
\tilde{M}_Z^-)\).

The fourth version, \(\widehat{\op{HF}}\), corresponds to the homology of a
fundamental domain; more precisely, \(H_*(\tilde{M}_Z^-,
t\tilde{M}_Z^-)\). 
%\(t\) being the generator of deck transformation in
%the negative direction.

Using topological arguments very special to this specific
two-dimensional situation, Ozsvath-Szabo showed that these Floer
homologies depend on \(z\) only through the spin-c structure.
\end{example*}

\section{Equivariant aspects of Seiberg-Witten-Floer theory}
Here is a reformulation of the construction of Kronheimer-Mrowka.
This formulation comes in particularly handy for our discussion on
filtrated connected sum formulae, in \S8 below.  

\subsection{The algebraic \(S^1\)-bundle}
Let \((C, \partial _C)\) be a complex of \(R\)-modules, and \(U\) be a degree \(-2\)
chain map on \(C\). We form the following new complex
\begin{equation}
(S_U(C), \partial_S)=(C\otimes R[y], \, \partial_C\otimes\sigma+U
\otimes y),
\end{equation}
where \(\op{deg}(y)=1\), \(y^2=0\), and 
the homomorphism \(\sigma: R[y]\to R[y]\) is defined by
\[\sigma(a+by)=a-by\quad \text{for \(a, b\in R\)}.\]

We call \(S_U(C)\) an {\em algebraic \(S^1\)-bundle}, due to the
following observation:
\begin{lemma*}
If \(C=C_*(B)\) is the chain complex of a manifold \(B\), and \(U_*:
H_{*}(B)\to H_{*-2}(B)\) agrees with the cap product with a cohomology
class \(u\in H^2(B)\), then \(H_*(S_U(C))\) is the homology of the
\(S^1\)-bundle over \(B\) with Euler class \(u\).
\end{lemma*}
To see this, view \(S_U(C)\) as a double complex, and notice that the
\(E^2\) term of the associated spectral sequence agrees with the
\(E^2\) term of the Serre spectral sequence of the aforementioned
\(S^1\) bundle. 

\subsection{Jones's formulation of equivariant homologies}
In retrospect, Ozsvath-Szabo's definition of the first three versions
of Heegaard Floer homologies is reminiscent of Jones's formulation of
\(S^1\)-equivariant homologies, which we sketch below. 

Let \(T\) be a \(S^1\) space. There are three versions of equivariant
homologies (cohomologies) that fit into fundamental long exact
sequences
\begin{gather*}
\cdots\to H_{S^1}^n(T)\to \hat{H}^n_{S^1}(T)\to G_{S^1}^{n+2}(T)\to
H_{S^1}^{n+1}(T)\to \cdots\\
\cdots\to G^{S^1}_n(T)\to \hat{H}_n^{S^1}(T)\to H^{S^1}_{n-2}(T)\to
G^{S^1}_{n-1}(T)\to \cdots
\end{gather*}
where \(\hat{H}\) denotes the Tate version (localized version), and
\(G\) is the co-Borel version which is dual to the usual (Borel)
equivariant (co)homology over \(R[U]\). 
%(\(G_*^{S^1}\) dual to \(H^*_{S^1}\)).

Modeling on the Serre spectral sequence of the fibration
\begin{equation}\label{jones-fib}
T\to T\times_{S^1}ES^1\to BS^1,
\end{equation}
J. D. S. Jones \cite{jones} wrote down the three versions of equivariant
(co)homologies in the following alternative way. 

Let \(S\) be the singular chain module or singular cochain algebra
of \(T\). The \(S^1\) action on \(T\),
\[g: S^1\times T\to T,\] 
equips \(S\) with a natural degree 1 chain map, \(J\): 
\[
J(x)=\begin{cases}(-1)^{\op{deg}(x)}g_* \delta(z\otimes x) &\text{when \(S\) is the singular chain module},\\
(-1)^{\op{deg}(x)}g^*x/z &\text{when \(S\) is the singular cochain
  module},
\end{cases}
\]
where \(\delta: C_{*}(S^1)\otimes C_{*}(T)\to C_{*}(S^1\times T)\) is the
Eilenberg-Zilber product, and
\(z\) is the fundamental 1-cycle of \(S^1\). 

Writing 
\begin{gather*}
V^- :=R[u],\, 
V^\infty := R[u, u^{-1}], \, 
V^+ := R[u, u^{-1}]/u R[u],
\end{gather*}  
we define:
\begin{equation}\label{def:E}
E^\bullet(S):=(S\otimes V^\bullet, \, \partial_S\otimes 1+J\otimes u)
\quad \text{for \(\bullet=-, \infty, +\)}. 
\end{equation}
(Jones's notation for
\(E^-, E^\infty, E^+\) are \(U^-, U^{\wedge}, U^+\) for homologies,
and \(V^-, V^{\wedge}, V^+\) for cohomologies respectively). 
\begin{lemma*}[Jones]
In the above notation, 
\begin{gather*}
H_*(E^-(S))=H^{-*}_{S^1}(T) \quad\text{when \(S\) is the singular
cochain algebra;}\\
H_*(E^+(S))=H_{*}^{S^1}(T) \quad\text{when \(S\) is the singular chain module.}
\end{gather*}
Furthermore, the long exact sequence induced by the short exact
sequence
\[
0\to uE^-(S)\to E^\infty(S)\to E^+(S)\to 0
\]
is precisely the fundamental exact sequence for equivariant
homologies/cohomologies, depending on whether \(S\) is the chain or
cochain module.
\end{lemma*}

\subsection{$\op{HM}^{tot}$, and $\hat{\op{HM}}, \bar{\op{HM}},
  \check{\op{HM}}$ as equivariant homologies}
We now combine the discussion in \S 5.1, 5.2 into a re-interpretation
of Kronheimer-Mrowka's construction.

Recall that Kronheimer-Mrowka's
Floer complexes \(\hat{\op{CM}}\), \(\check{\op{CM}}\) are modeled on
Morse theory of the (real) blow-up of an \(S^1\)-space along its 
fixed-point-set. Thus, they are analogs of the chain complex of the
base space of a \(S^1\)-bundle. We may thus apply the algebraic
\(S^1\)-bundle construction in \S5.1 to these complexes.
Let 
\begin{gather*}
\op{SM}^-:=S_U(\hat{\op{CM}}); \quad 
\op{SM}^+:=S_U(\check{\op{CM}}).
\end{gather*}
These chain complexes come equipped with a degree 1 chain map, namely,
multiplication by the nilpotent variable \(y\). This is precisely the
\(J\)-map above under the interpretation of \(S_U\) as the chain
complex of an \(S^1\)-bundle.
One may now apply the constructions \(E^-, E^{\infty}, E^+\) to
either \(\op{SM}^+\) or \(\op{SM}^-\) to obtain a sequence of
equivariant homologies.
\begin{lemma*}
{\bf (a)} There is a degree 0 chain map \(S_U(j): \op{SM}^+\to \op{SM}^-\)
commuting with the \(J\)-map, which induces
isomorphism in homologies.

{\bf (b)} \(E^+(\op{SM}^+)\) has the same homology as \(\check{\op{CM}}\);
\(E^-(\op{SM}^-)\) has the same homology as \(\hat{\op{CM}}\).
\(H_*(E^\infty(\op{SM}^-))=\hat{\op{HM}}\otimes_{R[U]}R[U,
U^{-1}]\); similarly for the plus-check version.

{\bf (c)} \(E^+(\op{SM}^-)\) has the same homology as \(\check{\op{CM}}\);
\(E^-(\op{SM}^+)\) has the same homology as \(\hat{\op{CM}}\).
\end{lemma*}
To see part (a),  
recall the maps 
\[
j: \check{\op{CM}}_*\to \hat{\op{CM}}_*;\quad
p: \hat{\op{CM}}_*\to \bar{\op{CM}}_{*-1};\quad
i: \bar{\op{CM}}_*\to \check{\op{CM}}_*
\]
defined by Kronheimer-Mrowka, which induce the fundamental long exact
sequence of Seiberg-Witten-Floer homologies. We define 
\begin{gather*}S_U(j): S_U(\check{\op{CM}})_*\to S_U(\hat{\op{CM}})_*,\\
S_U(p): S_U(\hat{\op{CM}})_*\to S_U(\bar{\op{CM}})_{*-1},\\
S_U(i): S_U(\bar{\op{CM}})_*\to S_U(\check{\op{CM}})_*\end{gather*}
to be the natural generalization of the above maps \(j, p, i\). For
example, in terms of the notation of \cite{KM:floer}, 
\[
S_U(j)=\left(\begin{array}{cc}1 &0\\ 0& \bar\partial_u^s+\bar{U}_u^s y\end{array}\right).
\]
A straightforward computation shows that these are indeed chain maps, and
they induce a long exact sequence in homologies
\begin{equation}\label{MC-seq}
\xymatrix{
&H_*(S_U(\bar{\op{CM}}))\ar@{->}[rd]^{S_U(i)_*}\\
H_*(S_U(\hat{\op{CM}}))\ar@{->}[ur]^{S_U(p)_*} &&H_{*}(S_U(\check{\op{CM}}))\ar@{->}[ll]^{S_U(j)_*}}
\end{equation}
%\begin{equation}
%\cdots \to H_{*}(S_U(\check{\op{CM}}))\stackrel{S_U(j)_*}{\to}
%H_*(S_U(\hat{\op{CM}}))
%\stackrel{S_U(p)_*}{\to} H_{*-1}(S_U(\bar{\op{CM}}))\stackrel{S_U(i)_*}{\to}
%H_{*-1}(S_U(\check{\op{CM}}))\to \cdots
%\end{equation}
Furthermore, \(S_U(j)\), \(S_U(p)\), \(S_U(i)\) are all
\(J\)-preserving (i.e. commuting with \(J\)), 
since in all three cases, the \(J\)-map is simply
multiplication by \(y\).

On the other hand, it is easy to see from a spectral sequence
calculation that
\begin{equation}\label{MC-bar}
H_*(S_U(\bar{\op{CM}}))=0.
\end{equation}

(\ref{MC-seq}) and (\ref{MC-bar}) together imply that \(S_U(j)_*\) is an isomorphism.

Part (b) is also consequence of simple spectral sequence
calculations. 

Part (c) is in fact a special case of \cite{jones} Lemma 5.2, which states
that:
\medbreak
If there is a \(J\)-preserving chain map from \(S_1\) to \(S_2\),
inducing an isomorphism from \(H_*(S_1)\) to \(H_*(S_2)\), then
the induced map from \(H_*(E^\bullet(S_1))\) to \(H_*(E^\bullet (S_2))\) is
also an isomorphism, for \(\bullet=-,\infty, +\).
\medbreak
This is a simple consequence of the observation that there is a
spectral sequence to compute the homology of \(E^\bullet(S)\) from the
homology of \(S\) (analog of the Serre
spectral sequence of the fibration (\ref{jones-fib})), which is natural with respect
to \(J\)-preserving chain maps. 

Applying the above statement to the \(J\)-preserving map \(S_U(j)\)
yields part (c). 

\bigbreak 
Because of part (a) of the above lemma, we may now define 
\begin{equation}
\op{HM}^{tot}:=H_*(S_U(\check{\op{CM}}))=H_*(S_U(\hat{\op{CM}})).
\end{equation}
Parts (b) and (c) of the previous lemma, together with the next lemma,
imply the following corollary:
\begin{cor*}
The fundamental exact sequence of the
Seiberg-Witten-Floer homologies of Kronheimer-Mrowka
agrees with the fundamental exact sequence of equivariant
homologies: 
\begin{equation}\label{eq:KM}\small
\xymatrix{
\llap{\(\cdots       \)}H_*(uE^-(\op{SM}^\pm))\ar@{->}[r]\ar@{=}[d]
&H_*(E^\infty(\op{SM}^\pm))\ar@{->}[r]\ar@{=}[d]
&H_*(E^+(\op{SM}^\pm)) \ar@{->}[r]\ar@{=}[d]
&H_{*-1}(uE^-(\op{SM}^\pm))\rlap{\(\cdots\)}\ar@{=}[d]\\
\llap{\(\cdots       \)}\hat{\op{HM}}_{*+1}\ar@{->}[r]
&\bar{\op{HM}}_*\ar@{->}[r]
&\check{\op{HM}}_*\ar@{->}[r]
&\hat{\op{HM}}_*\rlap{\(\cdots\)}}
\end{equation}
\normalsize
%(The first and last isomorphisms have degree 1).
\end{cor*}
\begin{lemma*}[Localization Theorem] We have
\[\hat{HM}(Y, \frak{s}; [w])\otimes_{R[U]}R[U,
U^{-1}]=\bar{HM}(Y, \frak{s}; [w])\]
\end{lemma*}
This is the analog of the familiar localization theorem in equivariant
theory (see e.g. \cite{AB}). In this context, the proof relies on the
nilpotence of the \(U\) action on \(\op{CM}^o\), the submodule of
\(\hat{\op{CM}}\) generated by irreducible critical points.
%\pf\(U\) is nilpotent on \(CM^\circ\), and on
%\(CM^u=\Z(\op{Crit}_{\text{Jacobi Torus}})\otimes_\Z \Z[u]\), it is
%multiplication by \(u\) modulo nilpotent terms. 
%This is because we are considering the case where CSD is
%\(\R\)-valued, and with a flow \(\Delta CSD\leq 0\), \(=0\) only when
%it is a flow between different eigenspaces over the same point in the
%Jacobi Torus. (Notice that when the critical points are only cyclically
%graded, the constraint \(\Delta\op{gr}=-2\) pick out a unique lifting
%of a critical point \(\tilde{q}\) so that \(\langle \tilde{q},
%\partial \tilde{p}\rangle\neq0\) for fixed \(\tilde{p}\).) 
%Next, comparing the \(u-u\) component of \(\hat{CM}\) and
%\(\bar{CM}\), the difference is \(\bar{\partial}^s_u\partial
%^u_s\). But this is again nilpotent because of \(\partial^u_s\): any
%flow going out of JT will decrease CSD. 
%Lastly, one needs to see that \(U\) action is invertible on
%reducibles (hence \(\bar{CM}\) is really a \(\Z[U, U^{-1}]\)-module,
%and increases eigenvalues (eventually reaching sufficiently low u). 
%\epf

\section{$HMT$: Seiberg-Witten-Floer theory with Taubes's perturbations}

To bridge the gulf of differences between the two theories as pointed
out above, our approach is to introduce an
intermediate object: a third set of Floer homologies,
\(\op{HMT}^\bullet\), which also come in four flavors, \(\bullet=-,
\infty, +, \wedge\), and the first 3 versions fit into a fundamental
long exact sequence. 
\(\op{HMT}\) should probably reads ``Heegaard-Monopole-Taubes'',
meaning that it is a variant of the SWT Floer theory sketched in \S2,
whose definition also involves some Heegaard ingredients: the choice
of a Heegaard splitting of \(Y\), and a filtration associated to a
1-cycle \(\underline{\gamma}_z\). 
%This Floer theory exhibits characteristics of both
%Heegaard and Seiberg-Witten theories, and 

We shall show that the equivalence of \(\op{HMT}\) with either
theory is easier to establish, and 
Conjecture 1.1 is thus broken into two:
\begin{conj}
Let \((Y, {\frak s})\) and \([w]=2\pi c_1(\frak{s})\) be as in
Conjecture 1.1. Then there are isomorphisms of \(R[U]\)-modules:
\[\begin{split}
{\rm\bf (a)} &\, \bar{\op{HM}}(Y, {\frak s}; [w])=\op{HMT}^\infty(\underline{Y}, \underline{\frak s});\, 
\hat{\op{HM}}(Y, {\frak s}; [w])=\op{HMT}^-(\underline{Y}, \underline{\frak s}; \underline{\gamma}_z);\\
&\, \check{\op{HM}}(Y, {\frak s}; [w])=\op{HMT}^+(\underline{Y},
\underline{\frak s}; \underline{\gamma}_z); \, \op{HM}^{tot}(Y, {\frak s}; [w])=\widehat{\op{HMT}}(\underline{Y},
\underline{\frak s}; \underline{\gamma}_z).\\
{\rm\bf (b)} &\, 
\op{HF}^\infty(Y, {\frak s})=\op{HMT}^\infty(\underline{Y}, \underline{\frak s});\,
\op{HF}^{\pm}(Y, {\frak s})=\op{HMT}^\pm(\underline{Y},
\underline{\frak s}; \underline{\gamma}_z); \\
&\, \widehat{\op{HF}}(Y, {\frak s})=\widehat{\op{HMT}}(\underline{Y},
\underline{\frak s}; \underline{\gamma}_z).\\
\end{split}\]
Furthermore, these isomorphisms are all natural with respect to the
fundamental sequences of \(\op{HM}^\bullet\), \(\op{HMT}^\bullet\),
and \(\op{HF}^\bullet\).
\end{conj}
As explained before, the \(U\)-map in the Heegaard Floer theory acts
by deck transformation, while in Seiberg-Witten-Floer theory, it is
cap (cup) product with the generator of \(H^*(BS^1)\). 
For the intermediate \(\op{HMT}\) theory, there are two natural module
structures, one from deck transformation, and the other from
equivariant theory. It turns out that these two module structures are
identical, so \(U\) above denotes either action in this theory. 
 
The following triangle best illustrates the relation among the three
Floer theories.
\[
\xymatrix{
&\op{HMT}^\bullet \ar@{=}[rd]^{\text{\, \, Conjecture 6.1(a)}}\\
\op{HF}^\bullet \ar@{=}[ur]^{\text{Conjecture 6.1(b)\, \, }}
&&\op{HM}^\bullet\ar@{.}[ll]^{\text{Conjecture 1.1}}}
\]

\subsection{The setup: $\underline{Y}$, $\underline{\gamma}_z$, and
  $\underline{\frak s}$}
We now describe the construction of \(\op{HMT}^\bullet\).
%As the notation suggests, this is in fact a
%variant of Seiberg-Witten (monopole) homologies, and `T' stands for
%Taubes, meaning that its definition uses Seiberg-Witten equations with
%perturbations of Taubes's type. The idea is another twist of the line of
%thinking outlined in \S2.

Let \(f: Y\to \R\) be a Morse function adapted to a Heegaard decomposition
of \(Y\), as in \S1.3.
As explained before, the obstruction to making \(f\)
harmonic is the existence of the maximum and minimum 
of \(f\), namely \(p_-, p_+\).
%\[p_+\in
%H_+, \quad p_-\in H_-\]
%respectively. 
Attach a 1-handle to \(Y\) along
the two points \(p_+, p_-\), and call the resulting manifold
\(\underline{Y}\).
\[
\underline{Y}=Y'\cup_{\partial B_{p_+}\sqcup\, \partial B_{p_-}}\left([-1,
1]\times S^2\right),\quad \text{where \(Y':=Y\backslash (B_{p_+}\sqcup
B_{p_-})\)}.
\]
One may now extend \(f\) to a {\em circle-valued} Morse function 
\[
\underline{f}: \underline{Y}\to S^1,
\]
which has no extrema. It is then possible to choose a metric on
\(\underline{Y}\), making \(\underline{f}\) harmonic \cite{calabi}. 

We shall use the notation \(\Sigma_H:=\underline{f}^{-1}(0)\) to denote the
Heegaard surface, to distinguish it from other genus \(g\) surfaces
\(\Sigma\subset \underline{Y}\).
Given an interval \(I\subset S^1\), \(\underline{Y}_I\subset
\underline{Y}\) will denote \(\underline{f}^{-1}I\), and
\(\underline{Y}_{]a, b[}:=\underline{Y}_{S^1\backslash (a, b)}\).
 
Let \(z\in \Sigma_H\), \(\gamma_z\subset Y\) be the base point and
associated 1-chain defined in \S1.3.  
Let \(\underline{\gamma}_z\subset \underline{Y}\) be a 1-cycle through
\(z\) completing \(\gamma_z\): more precisely, we choose it such that:
\begin{description}
\item [(\(\Gamma1\))] There are \(\kappa\in S^1\), \(\delta\in \R^+\),
  such that \(\underline{Y}_{[\kappa-\delta, \kappa+\delta]}\) is
  contained in the added 1-handle, and \(\underline{\gamma}_z\cap 
  \underline{Y}_{]\kappa-\delta, \kappa+\delta[}\) is a gradient flow line of
\(\underline{f}\) through \(z\in \Sigma_H\);
\item   [(\(\Gamma2\))]\(\underline{f}\) is monotone
along \(\underline{\gamma}_z\);
\item   [(\(\Gamma3\))]\(\underline{\gamma}_z\) avoids the ascending and descending
  manifolds from the critical points of \(\underline{f}\) in
  \(\underline{Y}\backslash \Sigma_H\).
\end{description}
The values \(\kappa, \delta\) will be chosen so that Claim 6.4 (3) below
holds.
In addition, we also require that 
\begin{equation}\label{z-class}
[*d\underline{f}]=\alpha \op{P.D.}[\underline{\gamma}_z]\in H^2 (\underline{Y}) \quad \text{for
  some constant \(\alpha>0\)}.
\end{equation}
This may be achieved either by fixing the metric, then 
choosing \(z\) so that the homology class
\([\underline{\gamma}_z]\) meets the requirement, or, 
fixing the class \([\underline{\gamma}_z]\), a closer examination of Calabi's argument in
\cite{calabi} shows that a metric can be found so that the class
\([*d\underline{f}]\) satisfies (\ref{z-class}) \cite{katz}.  
%and a simple path \(\gamma_z\) from \(p_+\)
%to \(p_-\) intersecting \(\Sigma_g\) at exactly one point, \(z\in
%\Sigma_g\). 
% \(\gamma_z\) extends over the 1-handle to form a 1-cycle
%\(\underline{\gamma}_z\) in \(\underline{Y}\). There is a harmonic \(S^1\)-valued function
%\(f:\underline{Y}\to \R/\Z\) adapted to the Heegaard decomposition, namely,
%\(f^{-1}[0, 1/2]\) contains \(g\) critical points in the interior, all of them of
%index 2, with boundary \(f^{-1}(1/2)\coprod
%(-f^{-1}(0))=S^2\coprod(-\Sigma_g)\). Similarly for the other half
%\(f^{-1}[1/2, 1]\). 

Let \(S_z\) be the boundary 3-sphere of a tubular neighborhood of
\(\gamma_z\). It splits \(\underline{Y}\) into a connected sum
\begin{equation}\label{z-split}
\underline{Y}=Y\#_{S_z}S^1\times S^2.
\end{equation}
Given a spin-c structure \(\frak{s}\) on \(Y\), let \(\underline{\frak
  s}\) be the spin-c structure on \(\underline{Y}\) defined by 
\begin{equation}\label{underline-s}
\underline{\frak s}={\frak s}\#_{S_z}{\frak s}_K,
\end{equation}
where \({\frak s}_K\) is the spin-c structure on \(S^1\times
S^2\) corresponding to the standard nowhere-vanishing vector field on
\(S^1\times S^2\), namely \(\nabla\tau\), \(\tau\) parameterizing
\(S^1\). 
Note that \({\frak s}_K\) is {\em not} the trivial spin-c structure. 
\[
c_1(\frak{s}_K)=2 \, \Omega_S,
\]
where 
%\(p_2: S^1\times S^2\to S^2\) is the projection, and
\(\Omega_S\) is the positive generator of \(H^2(S^1\times S^2)\).

\begin{rem*}
We use here and below the following orientation convention for
\(S^1\times S^2\) and its homologies: The parameterization
\(\tau\) orients \(S^1\); \(S^2\) is given the complex orientation,
and \(S^1\times S^2\) given the product orientation. \(H^1(S^1\times
S^2)\), \(H^2(S^1\times S^2)\) are oriented via their isomorphisms to
\(H^1(S^1)\), \(H^2(S^2)\) respectively. 
\end{rem*}

The Floer homologies \(\op{HMT}\) are constructed 
from the Seiberg-Witten-Floer
theory on \((\underline{Y}, \underline{\frak s})\), with perturbation of the form 
\begin{equation}\label{taubes-perp}
\omega=r*d\underline{f}+w,
\end{equation}
%(up to further small perturbations not changing the class of
%\(\op{PR}\)), 
where \(r\gg0\) is a constant, and \(w\) is a closed
2-form with cohomology class \([w_Y]\#_{S_z}[0]\), where \([w_Y]\)
satisfies (\ref{w-class}). 
For the rest of this article, a ``SWT'' Floer theory on
\(\underline{Y}\), or a ``SWT'' solution will refer to this particular
variant of SWT theory, unless otherwise specified. 

\subsection{Basic properties of \(\op{HMT}\) theory}

{\bf (1)} The decomposition (\ref{z-split}) splits \[H_1(\lc)=H^1(\underline{Y})=\Z\oplus
H^1(Y).\] 
In terms of this splitting, 
\[\begin{split}
\op{SF}&=c_1(\underline{\frak s})=2\oplus c_1({\frak s}), \\
\op{PR}&=2\pi(2\pi c_1(\underline{\frak s})-[w+r*df])=-C'\oplus 0,\quad
\text{for some constant \(C'>0\).}
\end{split}\]
Notice the complete agreement of these formulae with the formulae for
\(G, \op{SF}\), and \(\op{PR}\) in Heegaard Floer theory given in
\S3.2.  
\bigbreak

\noindent{\bf (2)} From the form of (\(SW_3\)), reducible critical points
exist only when 
\[
c_1(\frak{s})-[\omega]/(2\pi)=0.
\]
With our choice of \(\omega\), the left hand side is
\(-r/(2\pi)[*d\underline{f}]\), which is never zero as \(r>0\). 
Thus, all critical points (of the CSD functional) are irreducibles
(i.e. smooth points in \(\lc\)).

This again agrees with Heegaard Floer theory formally.  
\bigbreak

\noindent{\bf (3)} The above formula for \(\op{PR}\) indicates that it
always lies in the positive ray along the negative generator of \(H^1(S^1\times
S^2)\subset H^1(\underline{Y})\). 
Thus, the Morse-Novikov picture explained in \S3.1 leads to the
expectation that the Floer homology in this theory
(\(\op{HMT}^\infty\)) is an invariant,
namely, independent of \(r\), further exact perturbation of \(w\), and depends
on metric and \(\underline{f}\) only through the cohomology classes
\([d\underline{f}]\), \([*d\underline{f}]\). The compactness results
for Seiberg-Witten moduli spaces proven in \cite{KM:floer},
\cite{froy:cpt} confirm this expectation.

\subsection{Filtration by holonomy, and the definition of \(\op{HMT}^\bullet\)}
Because of the observation \S6.3 (2), we may set
\[
\op{CMT}^\infty(\underline{Y}, \underline{\frak
  s}):=\hat{\op{CM}}(\underline{Y}, \underline{\frak s}; [\omega])=\check{\op{CM}}(\underline{Y}, \underline{\frak s}; [\omega]),
\] 
with \(\omega\) given by (\ref{taubes-perp}). 
As we saw in \S6.3 (1), this Floer complex has exactly the same formal
properties as Ozsvath-Szabo's \(\op{CF}^\infty\). In particular, it
models on the Morse complex of the infinite-cyclic covering of \(\lc\)
determined by the homomorphism
\[
/[\underline{\gamma}_z]: H_1(\lc)=H^1(\underline{Y})\to \Z.
\]
We denote this infinite-cyclic covering by \(\hat{\lc}_{\underline{\gamma}_z}\).

To complete the analogy, we shall introduce a filtration on 
\(\op{CMT}^\infty\) in parallel to
the filtration by \(z\) in Heegaard Floer theory.

Notice that the holonomy of \(A^E\) 
along \(\underline{\gamma}_z\) defines a circle-valued function 
\[\op{hol}_{\underline{\gamma}_z}: \lc\to \R/2\pi\Z.\]
Let \(Z_{\underline{\gamma}_z}\) be the codimension 1 cycle in \(\lc\)
defined by
\begin{equation}
Z_{\underline{\gamma}_z}:=\op{hol}_{\gamma_z}^{-1}(c), \quad c\in
(\R/2\pi\Z )\backslash \{0\}.
\end{equation}
Say, let \(c=\pi\).
As explained in \S4.1, this 1-cocycle defines a local system on
\(\lc\). In fact, the Floer complex with this local coefficient is
precisely the Floer complex of the infinite-cyclic covering
\(\hat{\lc}_{\underline{\gamma}_z}\), since by (\ref{z-class}), 
\(Z_{\underline{\gamma}_z}\) represents the same cohomology class as
the above homomorphism \(/[\underline{\gamma}_z]\).

\begin{claim*}[1]
\(Z_{\underline{\gamma}_z}\) is semi-positive.
\end{claim*}
Thus, the construction in \S4.2 may be carried over to this context to
define the filtrated versions of Floer complexes
\[
\op{CMT}^-(\underline{Y},
\underline{\frak s}; \underline{\gamma}_z), \quad
\op{CMT}^+(\underline{Y},
\underline{\frak s}; \underline{\gamma}_z), \quad
\widehat{\op{CMT}}(\underline{Y},
\underline{\frak s}; \underline{\gamma}_z)
\]
corresponding respectively to the Floer complexes of the lower half of
\(\hat{\lc}_{\underline{\gamma}_z}\), of the pair, and of a
fundamental domain.  

It should now be clear that these filtrated complexes carry two
natural module structures over the polynomial ring of \(R\): we denote
by \(t\) the negative %(i.e. decreases the CSD functional) 
generator of deck transformation on \(\hat{\lc}_{\underline{\gamma}_z}\), and
reserve \(U\) for the degree \(-2\) chain map from \(S^1\)-equivariant
theory (defined from higher dimensional moduli spaces of flows).
These two module structures will be referred to as the \(R[t]\)- and 
\(R[U]\) module structures respectively, before they are finally
identified in \S8.1.
 
\begin{claim*}[2]
For \(r\gg1\), the corresponding homologies \(\op{HMT}^\bullet\) are independent
of the choice of level \(c\).
% exact perturbations to \(w\) of \(O(1)\), and \(r\).
\end{claim*}

The verification of these two claims requires a modification of Taubes
work in Part I of \cite{T}. 

Let \((\hat{A}^E, (\hat{\alpha}, \hat{\beta}))\) be a large \(r\) solution of the
SWT equations on \(\R\times \underline{Y}\). Recall that
according to Taubes's description of the local
behavior of SWT solutions, \(\hat{A}^E\) or \(A^E\) approximates 0
away from the corresponding Taubes curve or Taubes orbit.
%(possibly disconnected, with multiplicity) pseudo-holomorphic curve in
%\((\underline{Y}\backslash (d\underline{f})^{-1}(0))\times \R\).

%In particular, if \((\hat{A}^E, (\hat{\alpha}, \hat{\beta}))\) is 
%constant in \(s\) (which parametrizes the \(\R\)-factor), namely
%\((\hat{A}^E, (\hat{\alpha}, \hat{\beta}))=(A^E, (\alpha, \beta))\)
%for a solution to (\(\op{SW}_3\)), \((A^E, (\alpha, \beta))\), then 
%\(A^E\) approximates zero away from a set of gradient flow lines of 
%\(\underline{f}\) ending at \((d\underline{f})^{-1}(0)\). 

According to (\(\Gamma3\)), \(\underline{\gamma}_z\) avoids the
Taubes orbits; thus the holonomy map \(\op{hol}_{\underline{\gamma}_z}\)
always takes values near 0. 
This justifies our choice of the level \(c\).
Furthermore, note that two homologous codimension 1 cycles
\(Z\), \(Z'\) defines the same filtrated Floer complexes if they
together bound a region containing no critical points. 
The above observation on \(\op{hol}_{\underline{\gamma}_z}\) then
implies that different level sets of
\(\op{hol}_{\underline{\gamma}_z}\) defines the same filtration as
long as the level \(c\) stays away from 0.

To see Claim (1), we need the following interpretation of 
the holonomy filtration in terms of intersection numbers.
% of \(\R\times\underline{\gamma}_z\) with the Taubes curve. 

Let \(\tilde{\op{hol}}_{\underline{\gamma}_z}: \tilde{\lc}\to \R\) be
a lift of \(\op{hol}_{\underline{\gamma}_z}\), and 
\((\hat{A}^E, (\hat{\alpha}, \hat{\beta}))\) be a SWT flow line
between the SWT critical points \(c_-:=(A^E_-, (\alpha_-, \beta_-))\) and
\(c_+:=(A^E_+, (\alpha_+, \beta_+))\). Then 
\[
\tilde{\op{hol}}_{\underline{\gamma}_z}(c_+)-\tilde{\op{hol}}_{\underline{\gamma}_z}(c_-)=\int_{\R\times\underline{\gamma}_z}
F_{\hat{A}^E}\sim \# \, \left(\text{Taubes curve} \cap (\R\times \underline{\gamma}_z)\right).
\]
According to the next Claim and the condition (\(\Gamma1\)), 
the cylinder \(\R\times\underline{\gamma}_z\) intersects the Taubes
curve in a region where it is pseudo-holomorphic. 
Thus, the intersection number is non-negative. This implies the
semi-positivity of \(Z_{\underline{\gamma}_z}\). 
\begin{claim*}[3]
Let the spin-c structure be given by (\ref{underline-s}). Then there exist  
\(\kappa, \delta\) such that all Taubes orbits in \(\underline{Y}\) are contained in 
\(\underline{Y}_{]\kappa-\delta, \kappa+\delta[}\). Moreover, with
respect to appropriately chosen metric, all Taubes curves in \(\R\times \underline{Y}\)
are contained in 
\(\R\times \underline{Y}_{]\kappa-\delta, \kappa+\delta[}\).
\end{claim*}
An ``appropriately chosen metric'' can be one of the following: 
(i) when there is a holomorphic sphere of class \(*\times [S^2]\) 
in \[\R\times \underline{Y}_{(\kappa-\delta,
  \kappa+\delta)}=\R\times(\kappa-\delta,  \kappa+\delta)\times S^2
\] with respect to the associated almost complex
structure; (ii) when \(\underline{Y}_{(\kappa-\delta,
  \kappa+\delta)}\) is included in a long neck \([-L, L]\times S^2\subset
\underline{Y}\); or, according to the Atiyah-Floer picture explained
in \S7.2, 9.2 below, (iii) when \(\underline{Y}\) includes a long neck
along the Heegaard surface, \([-L, L]\times \Sigma\supset \Sigma_H\).
Notice that case (i) can always be arranged, since the sphere
\(\underline{f}^{-1}(\kappa)\) is a symplectic curve, and hence there exists
almost complex structures with respect to which it is
pseudo-holomorphic. 
%The condition on \(\underline{\frak s}\) implies
%that this pseudo-holomorphic curve avoids the added 1-handle, namely,
%it is supported in \(Y\times \R\subset \underline{Y}\times \R\) (see \S7.2). 
%In this region, \(\underline{\gamma}_z\times\R\) is
%pseudo-holomorphic by our assumption on
%\(\underline{\gamma}_z\). 
%First, recall that the 3-dimensional Seiberg-Witten-Taubes solutions
%on \(\underline{Y}\) corresponds to set of gradient flow lines of
%\(\underline{f}\) ending at critical points. Because \(\underline{f}\)
%is now circle-valued, these may include closed orbits. 
%However, with the choice of \(\underline{\frak s}\) in
%(\ref{underline-s}), these spin-c structure only span the subspace
%\[
%H^2(Y;\Z)\subset H^2(\underline{Y};\Z)
%\]
%in the space of all possible spin-c structures on \(\underline{Y}\).
%Closed orbits are excluded with this constraint on \(\underline{\frak
%  s}\).

To see that the Taubes orbits in \(\underline{Y}\) avoid \(\underline{Y}_{(\kappa-\delta,
  \kappa+\delta)}\), notice that the choice of \(\underline{\frak s}\) imposes the
homological constraint that the intersection number of the
Taubes orbit with any level surface \(S_\tau:=\underline{f}^{-1}(\tau )\),
\(\tau\in (\kappa-\delta, \kappa+\delta)\) is 0. Moreover, as the
constituent flow lines are
oriented by the gradient flow of \(\underline{f}\), they all intersect \(S_\tau\) positively. 
Thus, the Taubes orbit does not intersect any \(S_\tau\).

%Next, we examine the Seiberg-Witten-Taubes solutions on \(\R\times
%\underline{Y}\). As explained before, these correspond to Taubes
%curves, which has boundary along the lines \(\R\times
%\op{Crit}(\underline{f})\) and approaches the sets of flow lines just
%discussed as \(s\in \pm\infty\). (\(s\) parametrizes the \(\R\)
%factor). One can see similarly that a Taubes curve in the spin-c structure
%\(\underline{\frak s}\) avoids \(\R\times S\) in \(\R\times
%\underline{Y}\). This implies that \(\partial_{MT}\) also has the same
%geometric interpretation as the cylindrical situation considered in
%\S2.3. 

%To see that the Taubes curves in \(\R\times\underline{Y}\) avoids 
%\(\R\times \underline{Y}_{(\kappa-\delta, \kappa+\delta)}\), note that
The assertion on the Taubes curves can be seen as follows: In case
(iii), by reducing to the case of Taubes orbits; in case (ii), by
reducing to the case of a product complex structure on 
\(\R\times\underline{Y}_{(\kappa-\delta, \kappa+\delta)}\); in case
(i), by combining the 
the homological constraint from the spin-c
structure (which says that the intersection number of the Taubes curve
with any closed surface in \(\R\times \underline{Y}_{(\kappa-\delta,
  \kappa+\delta)}\) is zero), with the observation that in this case, there is a
2-parameter family of pseudo-holomorphic spheres covering  
\(\R\times\underline{Y}_{(\kappa-\delta, \kappa+\delta)}\).

%Suppose the metric is such that \(\R\times
%\underline{Y}_{(\kappa-\delta, \kappa+\delta)}\) contains a
%holomorphic sphere in the aforementioned case. (We omit discussion on
%other possible choices of metric). For example, this holds when the
%almost complex structure on \(\R\times
%\underline{Y}_{(\kappa-\delta, \kappa+\delta)}\) is the product
%complex structure. Alternatively, as \(S_\tau\) are symplectic curves,
%it is always possible to adjust the metric (hence the almost complex 
%structure) so that it is holomorphic. 

%According to index computation, such a holomorphic sphere occurs in a
%2-dimensional family, and their images cover \(\R\times
%\underline{Y}_{(\kappa-\delta, \kappa+\delta)}\), according to
%translation invariance. 

%Intersection number of the Taubes curve with any member of this family
%must be nonnegative, and hence according to the homological
%constraint, zero. 
%\footnote{perhaps omit this proof}
% indep of \(w\): dominanted by \(r*df\). indep of \(r\): taubes's
% conv thm
\begin{rem*}
{\bf (a)} Notice that the definition of \(\underline{\frak s}\), the splitting
\(H^1(\underline{Y})=H^1(Y)\oplus \Z\), and the filtration all depends
on the class \([\underline{\gamma}_z]\).
This is similar to the dependence of spin-c structure and filtration
on \(z\) in Heegaard Floer theory. 

{\bf (b)} It might appear that, by adding a 1-handle, we are forced from the
simpler (in the sense of \S 2.3) \(\R\)-valued Morse theory back to
the more complicated Morse-Novikov situation in \S2.2. 
Claim (3) above shows that, when the spin-c structure is chosen as 
(\ref{underline-s}), 
this is not the case, and the simple picture in
\S2.3 is retained. Notice that these spin-c structures only span the
subspace \(H^2(Y)\subset H^2(\underline{Y})\) in the space of all
spin-c structures. For general spin-c structures on
\(\underline{Y}\), one would indeed need the more complicated picture
in \S2.2.  
\end{rem*}

\subsection{A fundamental example: \((\underline{Y}, \underline{\frak
    s})=(S^1\times S^2, \frak{s}_K)\)}

The following is the simplest example of \(\op{HMT}^\bullet\), which
also plays an essential role in the proof of Conjecture 6.1 (a).

Let \(Y=S^3\), and \(\frak{s}\) be the unique spin-c structure on
\(S^3\). Then \(\underline{Y}=S^1\times S^2\), and \(\underline{\frak
  s}=\frak{s}_K\). Endow \(\underline{Y}\) with the product metric,
and notice that with this choice of \(\underline{\frak s}\), the line
bundle \(E\) is trivial.

There is an obvious \(S^1\)-valued harmonic function on
\(\underline{Y}\), namely
\[
\underline{f}=\tau, \quad \text{where \(\tau\) parameterizes \(S^1\).}
\]
The perturbation two form is now \(\omega=r*d\tau\), and there is an
obvious solution to the 3-dimensional Seiberg-Witten equations with
this perturbation:
\[
(A^E, (\alpha, \beta))=(0, (\sqrt{r}, 0)) \quad \text{with respect to
  a trivialization of \(E\)}.
\]
It is also not hard to see that this is the unique solution. 

Recalling that \(H_1(\lc)=H^1(\underline{Y})=\Z\), and \(\op{SF}=2\),
while \(\op{PR}=-C'<0\), 
this unique solution generates \(\op{CMT}^\infty\) as a \(R[t,
t^{-1}]\)-module.

Choose \(\underline{\gamma}_z=S^1\times \{z\}\) for a point \(z\in
S^2\). This is a gradient flow line of \(\underline{f}\). The map
\(\op{hol}_{\underline{f}}\) sends the unique solution above to 0.

Thus, we have
\[\begin{split}
\op{HMT}^- &=uR[u], \\
\op{HMT}^\infty &=R[u, u^{-1}], \\
\op{HMT}^+ &= R[u, u^{-1}]/(uR[u]), \\
\widehat{\op{HMT}} &=R,
\end{split}\]
where \(\deg u=-2\), and the deck transformation acts by
multiplication by \(u\).

Notice that this agrees both with \(\op{HF}^\bullet (S^3)\) and \(\op{HM}^\bullet
(S^3)\), confirming Conjecture 6.1 in this case. 

Next, we demonstrate that the \(U\)-map from equivariant theory agrees
with the deck transformation.%, namely, multiplication by \(u\).

Recall the following geometric interpretation of the \(U\)-map (see
e.g. \cite{D:floer} Lemma 7.6 for a Yang-Mills analog): 
Let \(p, q\) be two SWT critical points  with \(\ind
(p)-\ind (q)=2\), then
\[\langle q, Up\rangle=\deg \op{hol}_\nu, \]
where \(\op{hol}_\nu\) is the following holonomy map 
\[
\op{hol}_\nu: \cam (p, q)\to S^1,
\]
and \(\cam (p, q)\) is the 1-dimensional (reduced) moduli space
of SWT flow lines between \(p\) and \(q\), \(\op{hol}_\nu(c)\) is the holonomy
of \(\hat{A}^E\) along the path
\(\R\times \nu\), with respect to a chosen framing of \(E\) over
\(\nu\in \underline{Y}\), and \(c=(\hat{A}^E, (\hat{\alpha}, \hat{\beta}))\).
%is a SWT solution on \(\R\times\underline{Y}\) in the moduli space, 

In our situation, let \(p\) be the unique critical point described
above, and \(q=up\). The unreduced moduli space of flows between \(p\)
and \(q\) is \(\R\times S^1\), consisting of pull-backs
of vortex solutions of vortex number 1 on \(\R\times S^1\) to
\(\R\times S^1\times S^2\). The degree of \(\op{hol}_\nu\) can be computed from the
integral of \(F^{\hat{A}^E}\) over the cylinder \(\R\times S^1\times\op{pt}\),
namely the vortex number, 1. 

\section{Heuristic explanation of the equivalence}
Here is this author's best attempt (for now) 
at a conceptual explanation of the
somewhat mysterious relation (as highlighted in \S3.2) between the two 
theories. The actual proof of equivalence will not adhere to the heuristic
picture sketched below, as it is hard to make rigorous. However 
the picture does provide a useful guideline.

\subsection{From \(\op{HMT}\) to \(\op{HM}\)}

It has been gauge theorists' dream to understand Floer homologies as
the homologies of certain generalized spaces (pro-spectra?)
\cite{CJS}. Just as the Floer homology models on the Morse homology on
\(\tilde{\lc}\), this ``space'' models on the set of points in
\(\tilde{\lc}\) contained in finite-energy flow lines. This idea is 
difficult to realize, see however \cite{bauer, mano} for some recent
progress in this direction, in the Seiberg-Witten context.

Let's assume for the moment the existence of such objects:
suppose \(\op{HM}\) and \(\op{HMT}\) are \(S^1\)-equivariant
homologies of the ``\(S^1\)-spaces'' \(S_M\), \(S_{MT}\) respectively.
Due to this author's ignorance of homotopy theory, we shall regard
them as ordinary topological spaces.
The topological meaning of Conjecture 6.1 (a) then hinges on the
special properties of the ``space'' \(S_{MT}(S^1\times S^2,
\frak{s}_K)\), and the behavior of these spaces under connected sums
of 3-manifolds.  In fact, \(S_{MT}(S^1\times S^2,
\frak{s}_K)\) provides the mechanism that transforms the equivariant
Seiberg-Witten theory into the non-equivariant Heegaard Floer theory
of filtrated \(\Z\)-covers. 

From the computation in \S6.5, one expects \(S_{MT}(S^1\times S^2,
\frak{s}_K)\) to be, very roughly, of the (\(S^1\)-equivariant) homotopy type of 
the infinite dimensional Hopf sphere, ``\(S^{\infty}_{-\infty}\)'',
i.e. \(S^1\)-equivariant version of the
pro-spectrum \(\C P^\infty_{-\infty}\) (cf. e.g. Example 6.2 in \cite{CJS}).
Moreover, the filtration by deck transformation on \(\tilde{\lc}\)
induces a filtration on \(S_{MT}\): letting \(t^n\tilde{\lc}^-\)
denote the half of \(\tilde{\lc}\) below the hypersurface
\(t^n\tilde{Z}\), 
\[S_{MT}=\lim_{n\to -\infty}t^nS_{MT}^-, \] where 
\(t^nS_{MT}^-\) is defined from the flows in \(t^n\tilde{\lc}^-\). 
On the other hand, let \(S_{MT}^+=S_{MT}/S_{MT}^-\).

Both \(S_{MT}^\pm\) are homotopic to
\(S^\infty=ES^1\), and \(t^nS_{MT}^-\) is homotopic to \(t^{n+1}S_{MT}^-\) with 
a free \(S^1\)-cell \(e_S\) attached, 
\[e_S:=S^{2*+1}\backslash S^{2*-1}=\left\{v\, |\, v\in \C^{*+1}\backslash
\C^*, |v|=1\right\},\] 
where 
%\(t\) is a generator of the \(\Z\) action, and 
\(*\) is taken to the infinity in the limit. In addition, 
\begin{itemize}
\item Due to spectral flow, \(t^{n-1}S_{MT}^-\simeq \Sigma^\C t^nS_{MT}^-\),
  and \(S_{MT}=\Sigma^\C S_{MT}\).
\item The Euler class of the \(S^1\)
action is represented by a codimension 2 cycle
\(\mathcal{E}\), such that \(\mathcal{E}\cap t^nS_{MT}^-=t^{n+1}S_{MT}^-\). 
\end{itemize}
On the other hand, one expects the ``\(S^1\)-spaces'' of a connected sum of
3-manifolds to be a product 
\begin{gather*}
S_M(Y_1\# Y_2)\simeq S_M(Y_1)\times S_M(Y_2); 
\quad \text{in particular,}\\
S_{MT} (\underline{Y}, \underline{\frak s})\simeq S_M(Y, \frak{s};
[w])\times S_{MT}(S^1\times S^2, \frak{s}_K).
\end{gather*}
Notice that the \(S^1\) action on \(S_{MT}(S^1\times S^2,
\frak{s}_K)\) is free, implying that the \(S^1\)
action on the product \(S_{MT} (\underline{Y}, \underline{\frak
  s})\) is always free, even though the \(S^1\) action on
\(S_M(Y, \frak{s}; [w])\) has fixed points. Thus, the quotient is a
fiber product
\[
S_{MT} (\underline{Y}, \underline{\frak s})/S^1\simeq S_M(Y, \frak{s};
[w])\times_{S^1} S_{MT}(S^1\times S^2, \frak{s}_K), 
\]
which is smooth.
 
This explains why, in Seiberg-Witten-Taubes theory, all critical
points are irreducibles, and via Conjecture 6.1 (b),
why the Heegaard Floer theory is a
{\em non-equivariant} theory of a {\em smooth} \(\lc\).

Furthermore, the filtration on \(S_{MT}(S^1\times S^2,
\frak{s}_K)\) transfers to a filtration on \(S_{MT} (\underline{Y}, \underline{\frak
  s})\), and the diagonal \(S^1\) action on the product is determined by
the free \(S^1\) action on \(S_{MT}(S^1\times S^2,
\frak{s}_K)\). Thus, the agreement of deck transformation and \(U\) action
in \(\op{HMT}\) theory merely reflects the relation between Euler
class and the filtration via deck transformation on \(S_{MT}(S^1\times
S^2, \frak{s}_K)\), as described above. 
This in terms translates the equivariant aspects of 
Seiberg-Witten theory into the deck transformation in Heegaard Floer
theory. 

The hat versions of \(\op{HF}\) and \(\op{HMT}\) are both equivariant homologies
of flows in a fundamental domain. According to the above picture, They would
compute the homology of \[S_M(Y, \frak{s}; [w])\times_{S^1} e_S \simeq
S_M(Y , \frak{s}; [w]),\] 
namely the (non-equivariant) homology of the total \(S^1\)-space.

The \(+\) versions of \(\op{HF}\) and \(\op{HMT}\) are both
equivariant homologies of flows in
the upper half of a infinite cyclic cover. 
According to the above picture, they would
compute the homology of \[S_M(Y, \frak{s}; [w])\times_{S^1}
ES^1,\]
namely the equivariant homology of \(S_M(Y, \frak{s}; [w])\).

\subsection{From \(\op{HMT}\) to \(\op{HF}\)}

Because of Claim 6.4 (3), this is predicted by Taubes's
picture and a Seiberg-Witten analog of the Atiyah-Floer conjecture, 
similar to the 3-dimensional cylindrical SWT theory sketched in \S2.3. 
%We now examine this picture in the present context in more detail.
%Let \(\underline{Y}, \underline{\frak s}, \underline{f}\) be as
%described in \S6. 

%Finally, one may further reduce the geometric picture to 2-dimensions
%by stretching along the Heegaard surface as follows. 
Let 
\[
\underline{Y}(L)=\underline{Y}\backslash \Sigma _H\cup [-L/2,
L/2]\times \Sigma
\]
denote \(\underline{Y}\) endowed with a metric such that it contains a
cylinder of length \(L\) about the Heegaard surface \(\Sigma_H\), and
let \[Y(L)^\circ:=\underline{Y}(L)\Big\backslash
  \left(-\frac{L}{2+\epsilon}, \frac{L}{2+\epsilon}\right)\times \Sigma,\]
namely, a connected sum of the handlebodies \(H_\pm\) along
\(p_\pm\). 
Partition \(\R\times
\underline{Y}(L)\) into the union of \(\R\times [-L/2,
L/2]\times\Sigma\) and \(\R\times Y(L)^\circ\). We call the former the
{\em ``inside piece''}, the latter the {\em ``outside piece''}. 

The idea is that as \(L\to\infty\),
most of the ``energy'' is expected to reside on the inside piece.
Thus, the Taubes curve  in the
outside piece  would then approach a path of 
Taubes orbits
\begin{gather*}
\bigcup_{s\in \R}\{s\}\times \mathcal{O}_s, \quad \text{where \(\mathcal{O}_s\)
  is a Taubes orbit on \(Y^\circ\), and}\\
Y^\circ=Y(L)^\circ\cup\R_+\times \Sigma\cup \R_-\times \Sigma 
\end{gather*}
completes \(Y(L)\) into a 3-manifold with two cylindrical ends. 

Let \(\frak{s}^\circ\) denote the spin-c structure on \(Y^\circ\)
induced by \(\underline{\frak s}\).

A constituent flow line of a Taubes orbit in \(Y^\circ\) is specified
by its asymptotics at \(\tau\to\infty \) or \(-\infty\); in other
words, a point on the corresponding descending cycle on the limiting 
surface \(\Sigma_\infty\) at infinity, or ascending cycle on the
limiting surface \(\Sigma_{-\infty}\). With slight abuse of notation,
let \({\mathbb T}_\alpha\) be the product of descending cycles in
\(\sym^g\Sigma_+\), and \({\mathbb T}_\beta\) be the product of descending cycles in
\(\sym^g\Sigma_-\). Then a Taubes orbit in \(Y^\circ\) is specified by
a point in \({\mathbb T}_\alpha\times {\mathbb T}_\beta\). Namely,
there is a diffeomorphism
\[
\mathcal{O}_{Y^\circ}: {\mathbb T}_\alpha\times {\mathbb T}_\beta\to
\text{the space of Taubes orbits of spin-c structure
  \(\frak{s}^\circ\) in \(Y^\circ\).}
\]

%Moreover, expressing the Taubes curve as a union of 
%irreducible components, we
%claim that all the closed components are homotopic to level
%surfaces of \(\underline{f}\) in \(\underline{Y}\).
%To see this, notice that \(\underline{f}\) restricts to a
%\(S^1\)-valued harmonic function on any such component \(C_\nu\). 
%This means that \(\underline{f}|_{C_\nu}\) is either a constant, or a
%submersion. The latter possibility is excluded, since we have just
%seen that the Taubes curve never crosses \(\R\times S\).
%\footnote{the above holds only when the map \((f, s)\) from the nbhd in
%\(\R\times \underline{Y}\) to \(I\times \R\) is complex linear}

Next, we describe the more interesting inside piece. 

Note that the metric on \([-L/2, L/2]\times \Sigma \subset \underline{Y}(L)\)
is a product metric, and \(\underline{f}\) approximates a linear
function \(\tau+C_L\) as \(L\to\infty\). Thus, in the limit, the projection
\[
\text{inside piece} \to \R\times
[-L/2, L/2]
\]
is complex-linear. 
%with respect to the 
%the almost complex structure on \(\R\times
%\underline{Y}_{[-L/2, L/2]}\) associated with the metric and
%\(\omega\) and the standard complex structure on \(\R\times
%[-L/2, L/2]\). 
This implies that the Taubes curves in
the limit is a union of a multi-section of the \(\Sigma\)-bundle \(\R\times
[-L/2, L/2]\times \Sigma\), and a number of fibers. An index computation shows that
fibers do not appear as 
irreducible components of a Taubes curve in a reduced moduli space
of dimension \(\leq0\). On the other hand, we know that the
multi-section is a \(g\)-branched cover of \(\R\times
[-L/2, L/2]\), as in the \(s\to \pm \infty\) limit it is asymptotic to
a \[\text{Taubes orbit in \(\underline{Y}(L)\)} \cap [-L/2, L/2]\times \Sigma,
\]
which consists of \(g\) constituent flow lines. 
Furthermore, since this needs to match with the Taubes curve in the
outside piece, the multi-section  
restricts at the two boundary
components, \(\R\times\{-L/2, L/2\}\times\Sigma\) of the inside piece,
to the union of \(g\) paths on the \(g\) descending/ascending cycles. 

Thus, the Taubes curve of a large \(r\) SWT solution used for the
definition of the boundary map in the \(\op{HMT}\) Floer theory 
defines a holomorphic map:
\[
\mu_L: \R\times[-L/2, L/2]\to \sym^g\Sigma  \quad \text{with
  \(\mu_L(\cdot, -L/2)\in
  {\mathbb T}_\alpha\), \(\mu_L(\cdot,L/2)\in
  {\mathbb T}_\beta\)}.
\]
Composing with the conformal map \(\op{resc}_L^{-1}\), where
\[
\op{resc}_L: \R\times [-L/2, L/2]\to \R\times [0,1]: \, (s, \tau)\mapsto
(s/L, \tau/L+1/2),
\]
we obtain a holomorphic disk
\[
\mu: \R\times [0,1]\to \sym^g\Sigma, \quad \text{with \(\mu(\cdot, 0)\in
  {\mathbb T}_\alpha\), \(\mu(\cdot,1)\in
  {\mathbb T}_\beta\)},
\]
namely, a Heegaard flow line.

Conversely, a Heegaard flow line uniquely determines a limiting Taubes
curve. Thus, {\em if} an analog of Taubes's theorem (Theorem 2.1) holds in 
this context, this means \(\op{CMT}^\infty=\op{CF}^\infty\). 

The matching of filtration between the two Floer theories will not
require the full strength of Taubes's theorem, but only Part I of his
proof. We shall postpone its discussion to \S9.2.
\begin{rem*}
It follows that all irreducible components of a Taubes curve in this
context have multiplicity 1: All irreducible components with boundary has
multiplicity one, since we saw that the multiplicities at the boundary
are all 1. A closed irreducible component can only be a copy of \(\Sigma\), which
is of genus \(g>1\) according to our assumption. However, according to
Taubes, the multiplicity can be larger than 1 
only if \(g\leq 1\).
\end{rem*}
 
\section{Towards a real proof, part (a)}
This section and the next outline our plan for proving the two halves
of the conjecture, Conjecture 6.1 (a) and Conjecture 6.1 (b)
respectively.

Conjecture 6.1 (a) is now ``almost a theorem'', in the sense that 
details are still being written down \cite{L:model}, 
but there shouldn't be additional difficulties.

%\subsection{Conjecture 6.1 (a)}

Though the picture in \S7.1 is far from rigorous, it is
nevertheless possible to prove a connected sum formula for Floer {\em
  homologies} consistent with the fiber product picture, 
via a cobordism proof.
In the context of instanton Floer homologies, this is proven by Fukaya
\cite{fuk}, see also the exposition in \cite{D:floer}.

Because the Floer complexes \(\op{CM}^\bullet\) of Kronheimer-Mrowka
are built from real blown-ups of \(S^1\)-spaces along fixed-point-sets,
on which \(S^1\) acts {\em freely}, 
the connected sum formula in the context of Kronheimer-Mrowka
theory takes a cleaner form than the pre- blown-up version (as in 
\cite{fuk, D:floer}). The author learned of the following formulation
from Mrowka and Ozsvath. 
We only state the hat version, since it alone suffices for our purpose.
\begin{thm}[Connected sum formula]
Let \(Y_1, Y_2\) be closed, oriented 3-manifolds and \(\frak{s}_1\),
\(\frak{s}_2\) be spin-c structures on \(Y_1\), \(Y_2\)
respectively. For \(i=1, 2\), let \([w_i]\in H^2(Y_i;\R)\), and \(U_i:
\hat{\op{CM}}_*(Y_i, \frak{s}_i; [w_i])\to \hat{\op{CM}}_{*-2}(Y_i,
\frak{s}_i; [w_i])\) be the \(U\)-map defined in \cite{KM:floer}.
Then there is an isomorphism of \(R[U]\)-modules:
\[\hat{\op{HM}}_*(Y_1\# Y_2, \frak{s}_1\# \frak{s}_2;
[w_1]\#[w_2])=H_*\left(S_{U_1+U_2}(\hat{\op{CM}}(Y_1, \frak{s}_1; [w_1])\otimes_R
\hat{\op{CM}}(Y_2, \frak{s}_2; [w_2]))\right).
\]
\end{thm}
The \(R[U]\)-module structure of the right hand side of the above
isomorphism is given by the interpretation of the complex 
\[S_{U_1+U_2}(\hat{\op{CM}}(Y_1, \frak{s}_1; [w_1])\otimes_R
\hat{\op{CM}}(Y_2, \frak{s}_2; [w_2])\] 
as an ``algebraic fiber product''. 

Indeed, suppose \(C(B_1)\), \(C(B_2)\) are chain complexes of the base manifolds
of two \(S^1\)-bundles \(E_i\to B_i\) for \(i=1, 2\), and 
\(U_i: C_*(B_i)\to C_{*-2}(B_i)\) are
chain maps such that they induces maps on homologies that agree with 
cap products with the respective Euler class. Then 
\(C(B_1)\otimes C(B_2)\) is the chain complex of the base manifold of
\(E_1\times_{S^1} E_2\), and according to Lemma 5.1, 
\(H_*(S_{U_1+U_2}(C(B_1)\otimes C(B_2)))\) computes the homology of the
fiber product \(E_1\times_{S^1} E_2\).

The chain maps \[
U_1, -U_2 \, : S_{U_1+U_2}(C(B_1)\otimes C(B_2))\to S_{U_1+U_2}(C(B_1)\otimes C(B_2))
\]
are chain homotopic, and induce a degree \(-2\) map on homologies
agreeing with cap product with the Euler class of the \(S^1\)-bundle
\[
E_1\times E_2\to E_1\times_{S^1} E_2.
\]
%Let \([w]=[w_1]\#[w_2]\). In this case, there is a class \([z]\in
%H^2(X;\R)\) such that \([z]|_{Y_1\# Y_2}\), \([z]|_{Y_1}\),
%\([z]|_{Y_2}\) are \([w], [w_1], [w_2]\) respectively.  (\(X\) denotes
%the cobordism from \(Y_1\# Y_2\) to \(Y_1\coprod Y_2\)).
%Notice that the ``algebraic fiber product' \(MC\) is still a \(\Z[U]\) module: writing
%\(MC=C_1\otimes C_2\otimes \Z[y]\), \(\op{deg}(y)=1\), \(y^2=0\),
%\(\partial=(\partial_1+\partial_2)\sigma+(U_1+U_2)y\), then
%\(U=\pm U_1=\mp U_2+[\partial, \partial_y]\), where \(\sigma\) is a ``signature operator''
%\(\sigma(y)=-1\), \(\sigma(1)=1\). \footnote{check \(\pm\). change
%  notation: \(MC\to S\)}

The proof of Conjecture 6.1 (a) requires a filtrated version of the
above connected sum theorem.
Let 
\begin{gather*}
U_Y: \hat{\op{CM}}_*(Y, \frak{s}; [w])\to \hat{\op{CM}}_{*-2}(Y, \frak{s};
[w]),\\
U_{S^1\times S^2}: \op{CMT}^\bullet_*(S^1\times S^2, \frak{s}_K;
S^1\times \op{pt})\to \op{CMT}^\bullet_{*-2}(S^1\times S^2, \frak{s}_K; S^1\times \op{pt})
\end{gather*}
be respective the \(U\)-maps for \(Y\) and \(S^1\times S^2\).
For \(\bullet=-, \infty, +, \wedge\), denote by
\[
S_\otimes^\bullet (Y, \frak{s}):= S_{U_Y+U_{S^1\times S^2}}\left(\hat{\op{CM}}(Y, \frak{s}; [w])\otimes_R
\op{CMT}^\bullet(S^1\times S^2, \frak{s}_K; S^1\times \op{pt})\right).
\]
The first three of these complexes 
are endowed with both \(R[U]\)- and \(R[t]\)- module
structures, the former via its definition as an algebraic fiber
product; the latter via the \(R[t]\)-module structure of
\(\op{CMT}^\bullet(S^1\times S^2, \frak{s}_K; S^1\times \op{pt})\).
Furthermore, the short exact sequence
\[
0\to \op{CMT}^-(S^1\times S^2)\to
\op{CMT}^\infty(S^1\times S^2)\to 
\op{CMT}^+(S^1\times S^2)\to 0
\]
induces the short exact sequence
\[
0\to S^-_\otimes(Y, \frak{s})\to S^\infty_\otimes(Y, \frak{s})\to
S^+_\otimes (Y, \frak{s})\to 0,
\]
and hence a long exact sequence relating the homologies of \(S^\bullet_\otimes(Y, \frak{s})\).
\begin{thm}[Filtrated connected sum formula]
Let \(\bullet\) be any of the four: \(-, \infty, +\) or \(\wedge\),
and \(r\gg1\). Then there is an isomorphism 
\begin{gather*}
\op{HMT}^\bullet_*(\underline{Y}, \underline{\frak s};
\underline{\gamma}_z)
=H_*\left( S_\otimes^\bullet (Y, \frak{s})\right).\\
%\op{HMT}^\infty_*(\underline{Y}, [w+r*_3df])=H_*(MC(\hat{\op{CM}}(Y, [w])\otimes_\Z
%\op{CMT}^\infty(S^1\times S^2, [r*_3dt]); U_Y+U_{S^1\times S^2}).\\
%\op{HMT}^+_*(\underline{Y}, [w+r*_3df], \gamma_z)=H_*(MC(\hat{\op{CM}}(Y, [w])\otimes_\Z
%\op{CMT}^+(S^1\times S^2, [r*_3dt]); U_Y+U_{S^1\times S^2}).\\
%\widehat{\op{HMT}}_*(\underline{Y}, [w+r*_3df], \gamma_z)=H_*(MC(\hat{\op{CM}}(Y, [w])\otimes_\Z
%\widehat{\op{CMT}}(S^1\times S^2, [r*_3dt]); U_Y+U_{S^1\times S^2}).\\
\end{gather*}
both as \(R[U]\)- and as \(R[t]\)-modules.

Furthermore, these isomorphisms are natural in the sense that the
following diagram is commutative: 
\[
\small
\xymatrix{
\llap{\(\cdots       \)}\op{HMT}^-_{*}(\underline{Y}, \underline{\frak s};
\underline{\gamma}_z)\ar@{->}[r]\ar@{=}[d]
&\op{HMT}^\infty_*(\underline{Y}, \underline{\frak s}) \ar@{->}[r]\ar@{=}[d]
& \op{HMT}^+_* (\underline{Y}, \underline{\frak s};
\underline{\gamma}_z)\ar@{->}[r]\ar@{=}[d]
& \op{HMT}^-_{*-1}(\underline{Y}, \underline{\frak s};\underline{\gamma}_z)
\rlap{\(\cdots\)}\ar@{=}[d]\\
\llap{\(\cdots       \)} H_*(S_{\otimes}^-(Y, \frak{s}))\ar@{->}[r]
&H_*(S_{\otimes}^\infty(Y, \frak{s})) \ar@{->}[r]
& H_*(S_{\otimes}^+(Y, \frak{s}))\ar@{->}[r]
& H_{*-1}(S_{\otimes}^-(Y, \frak{s}))\rlap{\(\cdots\)}}
%\begin{CD}
%\cdots       \op{HMT}^-_{*}(\underline{Y}, \underline{\frak s};
%\underline{\gamma}_z)
%@>>> \op{HMT}^\infty_*(\underline{Y}, \underline{\frak s})
%@>>>
%\op{HMT}^+_* (\underline{Y}, \underline{\frak s};
%\underline{\gamma}_z)
%@>>> \op{HMT}^-_{*-1}(\underline{Y}, \underline{\frak s};
%\underline{\gamma}_z)
%\cdots\\
% @| @| @| @| \\
%\cdots       H_*(S_{\otimes}^-(Y, \frak{s}))@>>> H_*(S_{\otimes}^\infty(Y, \frak{s}))@>>>
%H_*(S_{\otimes}^+(Y, \frak{s})) @>>> H_{*-1}(S_{\otimes}^-(Y, \frak{s}))\cdots.\\
%\end{CD}
\]
\normalsize
%where \(\op{HMT}^\bullet\), \(S_\otimes^\bullet\) are short-hands for
%the two sides of the isomorphisms respectively; the first row is the
%fundamental exact sequence of \(\op{HMT}\), and the second row is the
%long exact sequence induced by the short exact sequence
%which in turn is induced, via a tensor product, from 
\end{thm}
\bigbreak
\noindent {\em Sketch of proof.}
To establish the isomorphism between the unfiltrated versions,
i.e. \(\bullet=\infty\), let \((Y_1, \frak{s}_1)=(Y, \frak{s})\) and
\(Y_2=(S^1\times S^2, \frak{s}_K)\), and \([w_1]=[w]\), \([w_2]=[r*_3
d\tau]\) in Theorem 8.1. 
The idea of the typical proof of 
connected sum formulae is to consider the pairs
\((V, \gamma_V)\), \((W, \gamma_W)\), where \(V\), \(W\) are
respectively the cobordisms giving the connected sum and 
its reverse operation.
\(\gamma_V\subset V\) and \(\gamma_W\subset W\) are paths between the \(Y_1\)-end and the
\(Y_2\)-end. In Morse-theoretic picture, there is a monotonically
decreasing Morse function on \(V\) with exactly one critical point, of
index 1, and \(\gamma_V\) is the descending manifold from this
critical point. \((W, \gamma_W)\) is obtained by reversing this Morse
function.

The pair \((V, \gamma_V)\) defines a chain map \(V_*\) from the Floer complex
of the connected sum to the algebraic fiber product. Conversely, \((W,
\gamma_W)\) defines a chain map \(W_*\) from the latter to the former.
Next, observe that the 
compositions of cobordisms \(W\cup_{Y_1\coprod Y_2} V\) and
\(V\cup_{Y_1\# Y_2} W\) are related to the product cobordisms,
respectively by a 1-surgery along the 1-cycle \(\gamma_V\cup
\gamma_W\), and collapsing the 2-sphere formed by the 2-dimensional
descending and ascending manifolds from the two critical points in 
\(V\cup_{Y_1\# Y_2} W\). 
One may then show that the compositions \(W_*\circ V_*\) and \(V_*\circ
W_*\) are chain homotopic to the identity, by proving 
some surgery formulae.

To prove the filtrated versions, 
one needs in addition to show that \(V_*\), \(W_*\)
preserve the filtration. The proof is a variant of the proof of the
semi-positivity of \(Z_{\underline{\gamma}_z}\).
%\footnote{the par below can be omitted} 
%In our context, there
%is a harmonic 2-form \(\varpi \) on the complete manifolds \(V\) and
%\(W\), which is asymptotic to \(*_3d\underline{f}\),
%\(*_3d\tau\), and \(0\) at the \(\underline{Y}\), \(S^1\times
%S^2\), and \(Y\) ends respectively. The chain maps \(V_*\), \(W_*\)
%are defined by counting solutions of Seiberg-Witten equations on \(V\)
%and \(W\), perturbed by \(\op{SD} (r\varpi+w_C)\), where \(w_C\) is a
%closed 2-form asymptotic to \(w\), \(0\), \(w_Y\) at the same three ends
%respectively.
%There is a cylinder \(\hat{\gamma}\) in \(V\) or \(W\) ending at the
%1-cycles \(\underline{\gamma}_z\) in the \(\underline{Y}\) and
%\(S^1\times S^2\) ends respectively. It is easy to choose metric so
%that \(\hat{\gamma}\) is pseudo-holomorphic except on \(\R_-\times
%\underline{Y}_{(\kappa-\delta, \kappa+\delta)}\) in the
%\(\underline{Y}\) end. The Taubes curve in this context is
%well-defined in a tubular neighborhood of \(\hat{\gamma}\), and it
%avoids a neighborhood of the non-holomorphic part of \(\hat{\gamma}\).

\begin{rem*}
The cobordism proof sketched above 
works when one of \(\op{CM}(Y_1)\), \(\op{CM}(Y_2)\) involves only
irreducibles. In general, \(W_*\) is not well-defined, due to compactness
problems. Instead, the plan of Mrowka et al. is to use naturality of
\(V_*\) with respect to the surgery exact sequences. 
\end{rem*}
\bigbreak

As we showed in \S6.5 that the \(R[t]\)- and \(R[U]\)-module
structures of \(\op{CMT}^\bullet (S^1\times S^2, \frak{s}_K; S^1\times
\op{pt})\) agree, this implies that the \(R[t]\)- and \(R[U]\)-module
structures of \(S_\otimes^\bullet (Y, \frak{s})\) also agree, and
hence, as an immediate corollary:
\begin{cor}
The \(U\)-action on \(\op{HMT}^\bullet_{*}(\underline{Y}, \underline{\frak s};
\underline{\gamma}_z)\) agrees with deck transformation. 
\end{cor}
This corollary may also be obtained via a geometric
argument generalizing the computation for the \(S^1\times S^2\) case in \S6.5.

As another immediate corollary, we have
\begin{cor}
Conjecture 6.1 (a) holds. 
\end{cor}
To see this, recall that with the standard choice of \(\underline{f}\)
and \(\underline{\gamma}_z\) for \(S^1\times S^2\) described in \S
6.5, 
\[
\op{CMT}^\bullet (S^1\times S^2, \frak{s}_K; S^1\times
\op{pt})=V^\bullet,
\]
where \(V^\bullet\) is as in (\ref{def:E}), and \(U_{S^1\times S^2}\)
is simply multiplication by \(u\).

Combining this with the filtrated connected sum formulae, the simple
fact that
\[
E^\bullet S_U(C)=S_{U+u}(C\otimes_R V^\bullet ) \quad \text{for a
  complex of \(R[U]\)-modules \(C\),}
%E^\infty S_U(C)=S_{U+x}(C\otimes \Z[x, x^{-1}]), \quad
%E^+S_U(C)=S_{U+x}(C\otimes \Z[x, x^{-1}]/x\Z[x]).
\]
and Corollary 5.3, Conjecture 6.1 (a) follows immediately.

\section{Towards a real proof, part (b)}
Progress towards the harder half of the conjecture, 
Conjecture 6.1 (b), is still in a very preliminary stage. 
We shall however present some partial results and ideas in this direction.

\subsection{Structure of moduli spaces in 3-dimensional
  Seiberg-Witten-Taubes theory}
 
%Recall the definition of \((\underline{Y}, \underline{\frak s})\) and \(\underline{f}\), and
%let \(\Sigma=\underline{f}^{-1}(0)\) be a Heegaard surface of genus
%\(g\). Notice that
%\[
%\partial(\underline{Y}\backslash \Sigma)=-\Sigma \coprod \Sigma.
%\]
%Attach two infinite cylinders at the above two boundary components to
%form a complete manifold, 
Let \(Y^\circ\), \(\frak{s}^\circ\) be as in \S7.2.
A metric on \(Y^\circ\) is called {\em cylindrical} if it
restricts to a fixed product metric on the attached
cylinders. Corresponding to each cylindrical metric, there is a
harmonic function \[f^\circ: Y^\circ\to \R, \quad 
\text{with the asymptotic condition (\ref{asym-tau})
  as \(\tau\to \pm\infty\),}\]
%where \(\tau\) parameterizes the \(\R\) factor of the cylindrical ends
%\(\R_\pm\times \Sigma\). 
%This function \(f^\circ\) is unique up to a
%constant, 
Furthermore, when \(L\) is large, the restrictions
\(f^\circ |_{Y(L)^\circ}\), \(\underline{f} |_{Y(L)^\circ}\)
approximate each other.

Let \(\cam_{Swt_3}(Y^\circ, \frak{s}^\circ; r, \eta)\) denote the
moduli spaces of the 3 dimensional Seiberg-Witten equations on the
cylindrical manifold \((Y^\circ, \frak{s}^\circ)\), with perturbation 
\[
\omega=r*df^\circ+\eta.
\]
This is a special case of the 3-d cylindrical SWT theory sketched in \S2.3.
Solutions to this perturbed Seiberg-Witten equations approach
``pull-backs'' of vortex solutions on \(E\to \Sigma\) exponentially as
\(\tau\to \pm \infty\) (see \cite{L:moduli,
  L1} for a definition). Thus, there are well-defined {\em end-point
  maps} from this moduli space to the moduli space of vortices on
\(E\to \Sigma\):
\[
\partial_\pm: \cam_{Swt_3}(Y^\circ, \frak{s}^\circ; r, \eta)\to \sym^g\Sigma.
\]
(Recall that \(c_1(\frak{s}^\circ)[\Sigma]=2\). This implies, via
  (\ref{S-split}), that \(\deg E=g\)).
\begin{prop*}
For generic pair of cylindrical metric and exact 2-form \(\eta\) of
\(O(1)\), the moduli space \(\cam_{Swt_3}(Y^\circ, \frak{s}^\circ; r,
\eta)\) is a compact, orientable manifold of dimension \(2g\).
%\[
%\cam_{Swt_3}(Y^\circ, \frak{s}^\circ; r, \eta)=\Big\coprod_{n} \cam_{Swt_3}^{2n}(Y^\circ, \frak{s}^\circ; r, \eta)
%\]
Furthermore, in the limit \(r\to \infty\), 
\[
\partial_-\times \partial_+: \cam_{Swt_3}(Y^\circ, \frak{s}^\circ; r,
\eta)\to \sym^g\Sigma\times \sym^g\Sigma 
\]
is a smooth, proper submersion of degree 1 to 
\({\mathbb T}_\alpha\times {\mathbb T}_\beta\subset \sym^g\Sigma\times
\sym^g\Sigma \). 
%(\({\mathbb T}_\alpha,  {\mathbb T}_\beta\) are as in \S7.2). 
\end{prop*}
See \cite{L:moduli} for a more precise statement. An analogous
theorem on the structures of moduli spaces
for manifolds with Euclidean ends may be found in \cite{L1}.
\bigbreak

\noindent{\em Sketch of proof.} 
Consider any general 3-manifolds
with cylindrical ends \(H\), and let \(\cam_{Swt_3}(H)\) be the moduli
space of 3-d SWT solutions.  
A simplified version of the arguments
in \cite{L1} establishes the smoothness, compactness, and
orientability properties of these moduli spaces.
An adaptation of Taubes's convergence theorem in part I of \cite{T}
describes the image of the limiting end-point maps in terms of 
products of descending/ascending cycles. 

The last statement of the Proposition on the
degree of the limiting \(\partial_-\times \partial_+\) follows from
gluing theorems for SWT moduli spaces on cylindrical manifolds, and
known computation of the Seiberg-Witten invariants of closed
3-manifolds.

Details will appear in \cite{L:moduli}. 

\begin{rem*}
Taubes's picture (cf. \S7.2) leads one to expect the map
\(\partial_-\times\partial_+\) to be a diffeomorphism, but the above
weaker statement will hopefully suffice for our purpose.
To prove the the stronger statement predicted by Taubes's picture, one
would need a counterpart of part III of \cite{T} (\(\op{Gr}\to
\op{Sw}\)) in this
situation. This is not an easy task; see however \cite{L1} for some progress.
\end{rem*}

\subsection{Comparing \(\op{CMT}\) and \(\op{CF}\)}
This subsection will be of even sketchier character, since details for
this part are not existent at this moment.

\begin{conj*}
For all sufficiently large \(L\), there is an isomorphism of \(R[U]\) modules
\[
\op{CMT}^\bullet(\underline{Y}(L), \underline{\frak
  s};\underline{\gamma}_z)=\op{CF}^\bullet(Y, \frak{s}).
\]
\end{conj*} 
To prove this conjecture, notice first that the identification of the
chain groups follows readily from the description of
\(\cam_{Swt_3}(Y^\circ)\) in the previous Proposition, a simple
gluing theorem for moduli of 3-dimensional SWT solutions
mentioned already in \S9.1 above, and the computation
of \(\op{SF}\) and \(\op{PR}\) of both sides explained before. Thus,
it remains to identify the spaces of flow lines on both sides. 

\bigbreak

\noindent{\bf (i) From Heegaard flow lines to SWT flow lines:} 
What this part requires is
a gluing theorem of the following form. 
According to \cite{OS1}, the Heegaard Floer homologies are invariant
under isotopies of the \(\alpha\)- and \(\beta\)-cycles. We may thus
choose \({\mathbb T}_\alpha\), \({\mathbb T}_\beta\) to be as in
Proposition 9.1, namely, as products of descending/ascending cycles of
the two surfaces at infinity of \(Y^\circ\).

Given a Heegaard flow line%, namely a holomorphic disk
\[
\mu: \R\times[0,1]\to \sym^g\Sigma, \quad \text{with \(\mu(\cdot,0)\in
  {\mathbb T}_\alpha\), \(\mu(\cdot,1)\in
  {\mathbb T}_\beta\),}
\]
the goal of the gluing theorem is to give a gauge-equivalence class of
SWT solution on \(\R\times\underline{Y}(L)\), for large \(L\).

Composing \(\mu\) with the rescaling map \(\op{resc}_L\) introduced in
\S7.2, 
%\[
%[-L/2, L/2]\times \R\to [0,1]\times \R: \, (\tau, s)\mapsto
%(\tau/L+1/2, s/L),
%\]
we obtain a family of vortex solutions on \(\Sigma\) parameterized by
\((s, \tau)\):
\[\left\{ (a_{(s, \tau)}, \nu_{(s, \tau)})\, |\, 
%\text{\((a__{(\tau,s)}, \nu_{(\tau,s)})\) is a vortex solution on
%\(\Sigma\)}, 
(s, \tau)\in\R\times [-L/2, L/2]\right\}.\] This defines a Seiberg-Witten configuration
on the inside piece, \(\R\times[-L/2, L/2]\times \Sigma\), such that 
\[
(\hat{A}^E, (\hat{\alpha}, \hat{\beta}))|_{\{(\tau, s)\}\times \Sigma}= (a_{(\tau,s)},
(\nu_{(\tau,s)}, 0)).
\]
When \(L\) is large, \(\partial_s(\hat{A}^E, (\hat{\alpha}, \hat{\beta}))\) and
\(\partial_\tau (\hat{A}^E, (\hat{\alpha}, \hat{\beta}))\) for such configurations are
small, and thus this gives an approximate SWT solution on the inside piece.

On the other hand, the restriction of \(\mu\) to the boundary \(\{0, 1\}\times \R\)
defines a map \(\partial \mu: \R\to {\mathbb T}_\alpha\times {\mathbb
  T}_\beta\), which lifts, under \(\partial_-\times \partial_+\), to paths in \(\cam_{Swt_3}(Y^\circ,
\frak{s}^\circ)\):
\[
c_s^j\in \cam_{Swt_3}(Y^\circ,
\frak{s}^\circ), \quad s\in \R, \, j=1, \ldots, k.
\] 
This in turn defines \(k\)
Seiberg-Witten configurations on \(\R\times Y^\circ\), so that
\[
(\hat{A}^E, (\hat{\alpha}, \hat{\beta}))|_{\{s\}\times Y^\circ}= c_{s/L}^j \quad
\text{for some \(j\).}
\]
These are again approximate SWT solutions, because
\(\partial_s(\hat{A}^E, (\hat{\alpha}, \hat{\beta}))\) is small due to the rescaling by \(L\). 

Noting that the above approximate solutions over the inside piece
and the outside piece
%\[\R\times Y(L)^\circ\subset \R\times Y^\circ, \quad \text{where
%  \(Y(L)^\circ:=\underline{Y}(L)\backslash
%  (-L/(2+\epsilon), L/(2+\epsilon))\times \Sigma\)}\] 
are close to each other on the overlap, 
a splicing construction then yields \(k\) approximate
SWT solutions over \(\R\times \underline{Y}(L)\).
A first task is then to perform an error estimate on these approximate
solutions showing:
\[
L\cdot \op{error}\to 0 \quad \text{as \(L\to \infty\).}
\]

To perturb the approximate solutions to true solutions, one needs
the deformation operator at these approximate solutions, \(D_j\),  to be
invertible, and that 
\[
L^{-1} \cdot D_j^{-1} \quad \text{is uniformly bounded.}
\]
This should follow from the fact that similar statements hold both
for the outside piece and the inside piece, as long as \(\mu\) is
nondegenerate: both the outside piece and the inside piece have 
a good Fredholm theory via a separation of variables argument 
(see \cite{T} part III and \cite{L1} for some example of this argument
in the Seiberg-Witten context). For the inside piece, this argument
reduces the Fredholm theory to the Fredholm theory of a 
Cauchy-Riemann operator on \(\R\times [-L/2, L/2]\) with (finite
dimensional) totally real boundary conditions. For the outside piece, 
the uniform boundedness follows from an eigenvalue estimate on
the compact manifold \(Y(L)^\circ\) with APS boundary conditions.
%\[
%\overline{Y(L)}^\circ:=\overline{\underline{Y}(L)\backslash \Sigma\times \{0\}}
%%  \cup \partial  (\underline{Y}(L)\backslash \Sigma\times \{0\}) 
%\quad  \text{with APS boundary conditions}.
%\]
(See e.g. \cite{CLM}).  

If the orientation works out, the last statement of Proposition 9.1 on
the degree of the end-point map would imply that 
the signed count of the
\(k\) SWT solutions corresponding to \(\mu\) equals
1.

\bigbreak
\noindent {\bf (ii) From SWT flow lines to Heegaard flow lines:}
This part requires a convergence theorem.
Very roughly, consider a sequence of
SWT solutions \(c_i\) on \(\R\times\underline{Y}(L_i)\) with
\(L_i\to \infty\).

Rescale the restrictions to the inside piece 
%\([-L_i/2,L_i/2]\times\R\times \Sigma\) 
by \(L_i^{-1}\) in the
\(\tau\) and \(s\) directions to get a sequence of configurations on
\(\R\times[0,1]\times \Sigma\). These configurations satisfy
equations that break up into two parts: the first being \(L_i^{-1}\)
times the Cauchy-Riemann equation, and the essential part of the 
second being the vortex equation on \(E\to\Sigma\). 
One thus expects that some adiabatic analysis akin to that in
\cite{DoS} to show that this gives, in the \(L\to \infty\) limit, a
pseudo-holomorphic disk in the moduli space of vortices, namely
\(\sym^g\Sigma\).
%, since \(\deg E=g\) by our choice of spin-c structure.

On the other hand, rescale the restrictions to the outside piece 
\(\R\times Y(L)^\circ\) in the \(s\) direction
by \(L_i^{-1}\). The \(\partial_s\) term in the
SWT equations is replaced by \(L_i^{-1}\partial_s\)
in the rescaled equations. Thus, one expects the \(L_i\to\infty\)
limit to yield a path of 3-dimensional SWT solutions,
namely, a path in \(\cam_{Swt_3}(Y^\circ, \frak{s}^\circ)\), and hence
a path in \({\mathbb T}_\alpha\times {\mathbb T}_\beta\) via the
end-point map \(\partial_-\times \partial_+\) described in Proposition
9.1. 

The above \(L\to\infty\) limits
for the outside and the inside piece have to match, as they both come
from SWT solutions that agree over the overlaps
\(\R\times (-L/2, -L/(2+\epsilon))\times\Sigma\) and
\(\R\times(L/(2+\epsilon), L/2)\times\Sigma\). 
Thus, the holomorphic map 
\[
\mu: \R\times [0, 1]\to \sym^g\Sigma
\]
arising from limits in the inside region must satisfy the boundary
condition 
\[
\mu(\cdot, 0)\in {\mathbb T}_\alpha, \quad \mu(\cdot, 1)\in {\mathbb T}_\beta.
\]
This is precisely a Heegaard flow line. 
\begin{rem*}
Notice that rescaling is also crucial for the desired compactness
results to be possible. While a variant of the usual compactness
results for Seiberg-Witten moduli spaces may very likely yield a
``local compactness theorem'' for the SWT solutions
in our context, ``global compactness'' for gauge-theoretic moduli
spaces over infinite cylinders \(\R\times M\) requires additional
decay estimates for the solutions in the ends \(s\to\pm \infty\). This
is typically achieved by an eigenvalue estimate for the relevant
deformation operator over \(M\): more precisely, the minimum of the
absolute values of the eigenvalues should be bounded away from zero. 
This will not hold when \(M\) itself contains an infinite cylinder.
Instead, our approach involves \(M\) with a long, but finite cylinder
of length \(L\). As \(L\to\infty\), the minimal absolute value of the
eigenvalues goes to 0. The effect of rescaling is that, instead of the
minimal absolute value, only a uniform lower bound on \[L \cdot (\text{
minimal absolute value of the eigenvalues})\]
is required for compactness.
\end{rem*}
\bigbreak 
\noindent {\bf (iii) Matching the filtrations:} 
We explained in \S6.4 that the holonomy filtration in \(\op{HMT}\)
theory is given by the intersection
number of the Taubes curve with \(\R\times \underline{\gamma}_z\).
Our task now is to compare this intersection number with the
intersection number giving the filtration in Heegaard Floer theory. 

For this purpose, we construct curves that approximate the large \(L\)
Taubes curve in the inside and outside pieces respectively.

Regarding a Heegaard flow line, 
\[
\mu: \R\times [0,1]\to \sym^g\Sigma, \quad \text{with \(\mu(\cdot, 0)\in
  {\mathbb T}_\alpha\), \(\mu(\cdot, 1)\in
  {\mathbb T}_\beta\)},
\]
as a multi-section of the
\(\Sigma\)-bundle \( \R\times [0,1]\times\Sigma\), its graph defines,
via composing with the holomorphic map 
\[
\R\times [-L/2,L/2]\times\Sigma\to \R\times[0,1]\times \Sigma: \, 
(s, \tau, w)\mapsto (s/L, \tau/L+1/2, w), 
\]
a curve in the inside piece \(C_\mu\subset \R\times [-L/2, L/2]\times
\Sigma\).

The restriction \(\partial \mu\), on the other hand, defines a surface
\(B_{\partial\mu}\)
with boundary at the zero locus of \(\op{SD}(*_3d\underline{f})\) in
the outside piece \(\R\times Y(L)^\circ\):
\[
B_{\partial \mu}\cap (\{s\}\times
Y(L)^\circ)=\left(\mathcal{O}_{Y^\circ}(\partial\mu(s/L)) \right) \cap Y(L)^\circ,
%\bigcup_{i=1}^g\left(\gamma^{-i}_{\mu(0, s/L)}\cup \gamma^{+i}_{\mu(1, s/L)}\right), 
\]
where \(\mathcal{O}_{Y^\circ}\) is the diffeomorphism introduced in \S7.2.
%where \(\gamma^{-i}_{(a_1, \ldots, a_g)}\) is the flow line from the
%\(i\)-th index 2 critical point of \(\underline{f}\) to the point
%\(a_i\) in the \(i\)-th descending cycle on \(\Sigma\).  \(\gamma^{+i}_{(b_1, \ldots, b_g)}\)
%is defined similarly for \((b_1, \ldots, b_g)\in {\mathbb T}_\beta\).

If the gluing construction in (i) works, 
\(C_\mu\) and \(B_{\partial\mu}\) approximates the Taubes curve
on the inside and outside piece respectively. However, via
(\(\Gamma3\)), 
the cylinder \(\R\times \underline{\gamma}_z\) does not intersect
\(B_{\partial\mu}\) in the outside piece. Thus, the intersection
number of the Taubes curve with \(\R\times \underline{\gamma}_z\) is
given by 
\[
\#\left((\R\times \gamma_z) \cap C_\mu\right)=\# \left(\{z\}\times
\sym^{g-1}\Sigma  \cap \op{image}(\mu)\right),
\]
which also gives the filtration in Heegaard Floer theory. 

\subsection{Final remarks} 

Two main technical components of this program are respectively 
analogs of the Atiyah-Floer conjecture, and Taubes's work. 
The following comments compare the variants needed in our context and
the ``original'' versions.
\bigbreak
 
\noindent {\bf (1)} In the Seiberg-Witten
context, the Atiyah-Floer conjecture should relate the Seiberg-Witten-Floer
homology with the Floer homology of Lagrangian intersections for 
\((P; L_1, L_2)\), where \(P\) is the moduli space of the
dimension-reduction of the Seiberg-Witten equations over the Heegaard
surface, and \(L_1, L_2\) are respectively images of the end-point
maps from the moduli spaces of 3-dimensional Seiberg-Witten solutions
over the two handlebodies \(H_+, H_-\).
 
An important difference, though, is that we work with a different
spin-c structure: the restriction of \(\underline{\frak s}\) on \(\underline{Y}\) minus the
1-handle does not extend across the two 3-balls to give a spin-c
structure on \(Y\). In particular, in our program,
\(c_1(\underline{\frak s})[\Sigma_H]=2\), while \(c_1({\frak
  s})[\Sigma_H]=0\) in the Atiyah-Floer conjecture.
We also work with a Taubes type perturbation which does not extend
across the two balls. 
%reduction to the Heegaard surface yields vortex solutions of vortex
%number \(g\), while in the Atiyah-Floer conjecture (for Seiberg-Witten
%theory), the vortex number would be \(g-1\). 

These differences make our program substantially easier than the 
Atiyah-Floer conjecture. As we saw, there are no irreducibles in
our program, and the ``Lagrangians'' and other moduli spaces in our story
are all smooth. 

In an earlier work \cite{OS2}, Ozsvath-Szabo defined a
``theta-invariant'' corresponding to the Euler characteristic of
Seiberg-Witten-Floer homology, modeling more directly with the
Atiyah-Floer conjecture. It remains interesting to understand the
corresponding Floer homology, which is yet to be defined.

\bigbreak

\noindent{\bf (2)} Though Taubes's picture has been the guiding 
principle for this work, our current plan does not involve running
the whole gamut of Taubes's fundamental work. As outline above, we
shall only need variants of Part I of
\cite{T} (\(\op{Sw}\to \op{Gr}\)), 
which is the more accessible part when \(\varpi^{-1}(0)\neq
\emptyset\), as pointed out in the beginning of \S2.2. 
Part II of \cite{T} defines \(\op{Gr}\), which in our context is
replaced by the work of Ozsvath-Szabo. Part III of \cite{T},
\(\op{Gr}\to \op{Sw}\), is in our plan replaced by the weaker result
Proposition 9.1 above. Part IV of \cite{T} compares the Kuranishi
structures of \(\op{Gr}\) and \(\op{Sw}\) theories, which is
complicated mainly due to curves with multiplicity \(>1\). 
As we saw that in our context, all irreducible
components of the Taubes curves have multiplicity 1, most of the
discussion in this part can therefore be avoided.  

\bigbreak
\noindent{\small {\bf Acknowledgments.} The author's debt to Taubes
  should be apparent from the frequent appearance of his name in this
  article. We are grateful to Kronheimer and Mrowka, 
for inspiring remarks, and for generously 
answering many questions on the details of their book in
preparation. We also benefited from helpful conversations with Bauer,
Manolescu, Ozsvath, and Szabo.}

\end{document}